\documentclass[reqno,12pt]{amsart}
\usepackage{amssymb,amsthm,amsmath,amsfonts,amstext,mathrsfs,
url,fancybox}
\usepackage{latexsym}
\usepackage{fancyhdr}
\usepackage{epsfig}
\usepackage[unicode]{hyperref}

\usepackage{pgf}
\usepackage{tikz}

\usepackage[utf8]{inputenc}
\usetikzlibrary{arrows,automata}
\usetikzlibrary{positioning}

\tikzset{
    state/.style={
           rectangle,
           rounded corners,
           draw=black, very thick,
           minimum height=2em,
           inner sep=2pt,
           text centered,
           },
}

\usepackage{verbatim}

\hypersetup{colorlinks=true,
linkcolor=blue,
citecolor=blue,
urlcolor=blue,
filecolor=blue,
bookmarksnumbered=true,
pdfstartview=FitH,
pdfhighlight=/N
}

\newtheorem{thm}{Theorem}
\newtheorem{cor}{Corollary}

\newtheorem{prop}{Proposition}
\newtheorem*{defin}{Definition}

\newtheorem*{prob12}{Problem 12}
\newtheorem*{prob13}{Problem 13}

\newcounter{noteno}
\newcounter{exam}

\newenvironment{Note}%
	{\refstepcounter{noteno}%
	\begin{small}
	\medbreak\par\noindent{{\bf Note~\thenoteno}.}}%
	{\hfill{$\Box$}\end{small}\par\medbreak}

	{\refstepcounter{exam}%
	\begin{small}
	\medbreak\par\noindent{{\bf Example~\theexam}.}}%
	{\hfill{$\Box$}\end{small}\par\medbreak}

\input xy
\xyoption{all}

\newcommand{\m}{\mathbf}

\newcommand{\bl}{\bullet}
\newcommand{\p}{\partial}
\def\d{\,{\rm{d}}}

\def\v{\,{\varsigma}}

\title[Projective superflows. II]
{Projective superflows. II. \\ $O(3)$ and the icosahedral group}
\author[G. Alkauskas]{Giedrius Alkauskas}
\address{Vilnius University, Department of Mathematics and Informatics, Naugarduko 24, LT-03225 Vilnius, Lithuania}
\email{giedrius.alkauskas@mif.vu.lt}

\begin{document}
\begin{abstract} Let $\m{x}\in\mathbb{R}^{n}$. For $\phi:\mathbb{R}^{n}\mapsto\mathbb{R}^{n}$ and $t\in\mathbb{R}$, we put $\phi^{t}=t^{-1}\phi(\m{x}t)$. \emph{A projective  flow} is a solution to the projective translation equation $\phi^{t+s}=\phi^{t}\circ\phi^{s}$, $t,s\in\mathbb{R}$. \emph{The projective superflow} is a projective flow with a rational vector field which, among projective flows with a given symmetry, is in a sense unique and optimal.\\
\indent In this second part we classify $3-$dimensional real superflows. Apart from the superflow $\phi_{\widehat{\mathbb{T}}}$ (with a group of symmetries being all symmetries of a tetrahedron) and the superflow $\phi_{\mathbb{O}}$ (with a group of symmetries being orientation preserving symmetries of an octahedron), both described in the first part of this study, here we investigate in detail the superflow $\phi_{\mathbb{I}}$ whose group of symmetries is the icosahedral group $\mathbb{I}$ of order $60$. This superflow is a flow on co-centric spheres, and is also solenoidal, like the first two. These three superflows is the full (up to linear conjugation) list of irreducible $3-$dimensional real projective superflows.\\
\indent We also find all reducible $3$-dimensional real superflows. There are two of them: one with group of symmetries being all symmetries of a $3$-prism (group of order $12$), and the second with a group of symmetries being all symmetries of a $4$-antiprism (group of order $16$). 
\end{abstract}
\pagestyle{fancy}
\fancyhead{}
\fancyhead[LE]{{\sc Projective superflows}}
\fancyhead[RO]{{\sc G. Alkauskas}}
\fancyhead[CE,CO]{\thepage}
\fancyfoot{}

\date{September 20, 2016}
\thanks{The research of the author was supported by the Research Council of Lithuania grant No. MIP-072/2015}

\maketitle
\section{Projective flows}
\subsection{Preliminaries}
\label{prelim}
We skip the introduction, since the necessary details are contained in the first part of our study. We only recall that \emph{the projective translation equation} was first introduced in \cite{alkauskas-t} and is the equation of the form

\begin{eqnarray}\setlength{\shadowsize}{2pt}\shadowbox{$\displaystyle{\quad
\frac{1}{t+s}\,\phi\big{(}\m{x}(t+s)\big{)}=\frac{1}{s}\,\phi\Big{(}\phi(\m{x}t)\frac{s}{t}\Big{)}},\quad t,s\in\mathbb{R}\text{ or }\mathbb{C}$.\quad}\label{funk}
\end{eqnarray}

A non-singular solution of this equation is called \emph{a projective flow}.  The \emph{non-singularity} means that a flow satisfies the boundary condition
\begin{eqnarray}
\lim\limits_{t\rightarrow 0}\frac{\phi(\m{x}t)}{t}=\m{x}.
\label{init}
\end{eqnarray}

A smooth projective flow $\phi=u\bl v\bl w$ is accompanied by its \emph{vector field}, which is found from
\begin{eqnarray}
\varpi(x,y,z)\bl\varrho(x,y,z)\bl\sigma(x,y,z)=\frac{\d}{\d t}\frac{\phi(xt,yt,zt)}{t}\Big{|}_{t=0}.
\label{vec}
\end{eqnarray}
Vector field of a projective flow is necessarily a triple of $2$-homogenic functions. If a flow is smooth, the functional equation (\ref{funk}) implies the PDE \cite{alkauskas}
\begin{eqnarray}
u_{x}(\varpi-x)+u_{y}(\varrho-y)+u_{z}(\sigma-z)=-u,\label{pde}
\end{eqnarray}
and the same PDE for $v$ and $w$, with the boundary conditions as given by (\ref{init}). We emphasize that apart from the first paper \cite{alkauskas-t} where topologic aspects of projective flows were investigated, assuming only continuity but not smoothness of a flow, all the rest papers \cite{alkauskas, alkauskas-un,alkauskas-ab,alkauskas-comm, alkauskas-super1} investigate rational vector fields. So this area of research is in the intersection of algebraic geometry (birational geometry, curves), differential geometry (local and global flows, differential systems), number theory (abelian integrals and functions, algebraic functions, number fields) and algebra (group representations).

\subsection{General setting for $3-$dimensional flows} In this subsection we briefly recall the method to tackle PDE (\ref{pde}) in dimension $3$ \cite{alkauskas-super1}.  \\
  
Let $\varpi\bl\varrho\bl\sigma$ be a triple of $2-$homogenic rational functions, giving rise to the projective flow $u\bl v\bl w$. With the PDE (\ref{pde}) and the boundary condition (\ref{init}), we associate an autonomous system of ODEs:
\begin{eqnarray}
\left\{\begin{array}{l}
p'=-\varpi(p,q,r),\\
q'=-\varrho(p,q,r),\\
r'=-\sigma(p,q,r).
\end{array}
\right.
\label{sys-in}
\end{eqnarray}
Let $u(x,y,z)$ be the solution to the PDE (\ref{pde}) with the boundary condition as in the first equality of (\ref{init}), that is, $\lim\limits_{t\rightarrow 0}\frac{u(xt,yt,zt)}{t}=x$. Similarly as in \cite{alkauskas-un,alkauskas-super1}, we find that the function $u$ satisfies
\begin{eqnarray}
u\Big{(}p(s)\v,q(s)\v,r(s)\v\Big{)}=p(s-\v)\v,\quad s,\v\in\mathbb{R}.
\label{3-var}
\end{eqnarray}
As in \cite{alkauskas-super1}, suppose that the system (\ref{sys-in}) possesses two independent homogeneous first integrals $\mathscr{W}(p,q,r)$ and $\mathscr{V}(p,q,r)$: the homogeneous function $\mathscr{W}$ satisfies
\begin{eqnarray}
\mathscr{W}_{x}\varpi+\mathscr{W}_{y}\varrho+\mathscr{W}_{z}\sigma=0,
\label{first-int}
\end{eqnarray}
analogously for $\mathscr{V}$. To fully solve the flow in explicit terms, together with the system (\ref{sys-in}), we require
\begin{eqnarray*}
\mathscr{W}(p,q,r)=1,\quad \mathscr{V}(p,q,r)=\xi,
\end{eqnarray*}
where $\xi$ is arbitrary, but fixed. For a function $u$ in $3$ variables, (\ref{3-var}) now involves three parameters $s,\v$ and $\xi$. Generically, this gives the analytic formula for $u$.
\subsection{Finite subgroups of $O(3)$}
\label{sub-finite}
In order to find all $3$-dimensional superflows, we need to know all possible symmetry groups. The Table \ref{table1} lists all finite subgroups $\Gamma$ of $O(3)$. This is a classical list, hence we provide only a very short explanation.
\begin{itemize}
\item[i)]\emph{Rotation groups}. These are finite subgroups of $SO(3)$.
\item[ii)]\emph{Direct products}. If $\Gamma$ contains the matrix $-I$, then $H=\Gamma\cap SO(3)$ is a subgroup of $\Gamma$ of index $2$, and is a rotation group. Since $-I$ commutes with all matrices, we get a direct product: $\Gamma=H\times\{I,-I\}$. 
\item[iii)]\emph{Mixed groups}. In this case $\Gamma\not\subset SO(3)$, but $-I\notin\Gamma$. Such groups $\Gamma$ are described by a pair of rotations groups $KH$, $K$ containing $H$ as a subgroup of half size. Then the element of $\Gamma$ are either elements of a smaller group $H$, or they are equal to $-I\times T$, where $T\in K\setminus H$.
\end{itemize} 

\begin{table}[h]
\begin{tabular}{|c | c|c|| c| c|| c| c|}
\hline
\multicolumn{7}{|c|}{\textbf{Finite subgroups of $O(3)$}}\\
\hline\hline
\multicolumn{3}{|c||}{Rotation groups} & \multicolumn{2}{|c||}{Direct products} & \multicolumn{2}{|c|}{Mixed groups} \\
\hline\hline
Name & Group & Order & Group & Order & Group & Order\\
\hline
Cyclic & $\mathbb{C}_{d}$ & $d$ & $\mathbb{C}_{d}\times\{I,-I\}$ & $2d$ & $\mathbb{C}_{2d}\mathbb{C}_{d}^{(4)}$ & $2d$\\
\hline
Dihedral & $\mathbb{D}_{d}^{(6)}$ & $2d$ & $\mathbb{D}_{d}\times\{I,-I\}$ & $4d$ & $\mathbb{D}_{d}\mathbb{C}_{d}^{(5)}$ & $2d$\\
\hline
Tetrahedral & $A_{4}\simeq\mathbb{T}$ & $12$ & $\mathbb{T}\times\{I,-I\}$ & $24$ & $S_{4}A_{4}\simeq\widehat{\mathbb{T}}^{(1)}$ & $24$\\
\hline
Octahedral & $S_{4}\simeq\mathbb{O}^{(2)}$ & $24$ &  $\mathbb{O}\times\{I,-I\}$ & $48$ & $\mathbb{D}_{2d}\mathbb{D}_{d}^{(7)}$  & $4d$\\
\hline
Icosahedral & $A_{5}\simeq\mathbb{I}^{(3)}$ & $60$ & $\mathbb{I}\times\{I,-I\}$ & $120$ & \multicolumn{2}{|c|}{-}\\
\hline
\hline
\end{tabular}
\caption{List of finite subgroups of $O(3)$. The superscript $\Gamma^{(x)}$, $x\in\{1,\ldots,7\}$ is an order of a superflow in Theorem \ref{thm1}, and for $x>3$ is just for convenience to label the group and show later that superflow for it exists or does not exist. We will use these in Subsection \ref{results}.}
\label{table1}
\end{table}
\subsection{Superflows}
We briefly recall the notions of the superflows introduced in \cite{alkauskas-un,alkauskas-super1}.

\begin{defin}
Let $n\in\mathbb{N}$, $n\geq 2$, and $\Gamma\hookrightarrow{\rm GL}(n,\mathbb{R})$ be an exact representation of a finite group, and we identify $\Gamma$ with the image. We call the flow $\phi(\m{x})$ \emph{the $\Gamma$-superflow}, if 
\begin{itemize}
\item[i)]there exists a vector field $\mathbf{Q}(\m{x})=Q_{1}\bl\cdots\bl Q_{n}\neq 0\bl\cdots\bl 0$ whose components are $2$-homogenic rational functions and which is exactly the vector field of the flow $\phi(\m{x})$, such that
\begin{eqnarray}
\gamma^{-1}\circ\mathbf{Q}\circ\gamma(\mathbf{x})=
\mathbf{Q}(\mathbf{x})\label{kappa}
\end{eqnarray}
 is satisfied for all $\gamma\in\Gamma$, and
\item[ii)] every other vector field $\mathbf{Q}'$ which satisfies (\ref{kappa}) for all $\gamma\in\Gamma$ is either a scalar multiple of $\mathbf{Q}$, or  its degree of a common denominator is higher than that of $\mathbf{Q}$. 
\end{itemize} 
The superflow is said to be \emph{reducible or irreducible}, if the representation  $\Gamma\hookrightarrow{\rm GL}(n,\mathbb{R})$ (considered as a complex representation) is reducible or, respectively, irreducible.
\end{defin} 
Thus, if $\phi$ is a superflow,  then it is uniquely defined up to conjugation with a linear map $\m{x}\mapsto t\m{x}$. This corresponds to multiplying all components of $\mathbf{Q}$ by $t$. Some words are needed to explain what (\ref{kappa}) means. Simply, we interpret $\m{x}$ as a vector-column, $\gamma$ - as a linear map $\mathbb{R}^{n}\mapsto{\mathbb{R}}^{n}$ acting on vectors-columns, and the vector field $\mathbf{Q}$ - as a rational $2$-homogeneous map $\mathbb{R}^{n}\mapsto{\mathbb{R}}^{n}$ acting again on vectors-columns. So, $\gamma^{-1}\circ\mathbf{Q}\circ\gamma$ is again a $2$-homogeneous map $\mathbb{R}^{n}\mapsto{\mathbb{R}}^{n}$, which must be identically equal to $\mathbf{Q}$ for all $\gamma\in\Gamma$. \\
 
 There exist reducible superflows: over $\mathbb{C}$, in dimension $2$ \cite{alkauskas-super3}, and in this paper we will show that reducible superflows over $\mathbb{R}$ exist in any dimension $n\geq 3$; see Proposition \ref{prop-ext}.\\

In a $2$-dimensional case, we found that for every $d\in\mathbb{N}$, there exist the superflow whose group of symmetries is the dihedral groups $\mathbb{D}_{2d+1}$. This list exhaust all $2-$dimensional superflows. 
\subsection{Results}
\label{results}
Let 
\begin{eqnarray*}
\phi=\frac{1+\sqrt{5}}{2}.
\end{eqnarray*}
The first main result of this part reads as follows. 
\begin{thm}
\label{thm1} There exist the following three irreducible projective superflows over $\mathbb{R}$ in dimension $3$.
\begin{itemize} 
\item[1)]The superflow $\phi_{\widehat{\mathbb{T}}}$ with the vector field
\begin{eqnarray*}
yz\frac{\p}{\p x}+xz\frac{\p}{\p y}+xy\frac{\p}{\p z},
\end{eqnarray*}
whose group of symmetries is the group of order $24$, the symmetries of a tetrahedron. The orbits are space curves $\{x^2-y^2=\mathrm{const.}, x^2-z^2=\mathrm{const.}\}$, and generically these are curves of genus $1$. The flow can be integrated explicitly in terms of Jacobi elliptic functions.
\item[2)]The superflow $\phi_{\mathbb{O}}$ with the vector field
\begin{eqnarray*}
\frac{y^3z-yz^3}{x^2+y^2+z^2}\frac{\p}{\p x}+\frac{z^3x-zx^3}{x^2+y^2+z^2}
\frac{\p}{\p y}+\frac{x^3y-xy^3}{x^2+y^2+z^2}\frac{\p}{\p z},
\end{eqnarray*}
whose group of symmetries is the group of order $24$, orientation preserving symmetries of an octahedron. The orbits are space curves $\{x^2+y^2+z^2=\mathrm{const.}, x^4+y^4+z^4=\mathrm{const.}\}$, and generically these are curves of arithmetic genus $9$. The flow can be integrated explicitly in terms of Weierstrass elliptic functions, but in a rather complicated way, via reduction of special hyper-elliptic functions of genus $2$ to elliptic ones.
\item[3)]The superflow $\phi_{\mathbb{I}}$ with the vector field
\begin{eqnarray*}
\varpi\frac{\p}{\p x}+\varrho\frac{\p}{\p y}+\sigma\frac{\p}{\p z},
\end{eqnarray*}
where\footnotesize 
\begin{eqnarray}
\varpi&=&\frac{(5-\sqrt{5})yz^5+(5+\sqrt{5})y^5z-20y^3z^3+(10+10\sqrt{5})x^2yz^3+
(10-10\sqrt{5})x^2y^3z-10x^4yz}{(x^2+y^2+z^2)^2},\nonumber\\
\varrho&=&\frac{(5-\sqrt{5})zx^5+(5+\sqrt{5})z^5x-20z^3x^3+(10+10\sqrt{5})y^2zx^3+(10-10\sqrt{5})y^2z^3x-10y^4zx}{(x^2+y^2+z^2)^2},\label{vec-field}\\
\sigma&=&\frac{(5-\sqrt{5})xy^5+(5+\sqrt{5})x^5y-20x^3y^3+(10+10\sqrt{5})z^2xy^3+(10-10\sqrt{5})z^2x^3y-10z^4xy}{(x^2+y^2+z^2)^2}\nonumber
\end{eqnarray}\normalsize 
whose group of symmetries is the group  of order $60$, orientation preserving symmetries of an icosahedron. The orbits are space curves $\{x^2+y^2+z^2=\mathrm{const.}, (x^2\phi^2-y^2)(y^2\phi^2-z^2)(z^2\phi^2-x^2)=\mathrm{const.}\}$, and generically these are curves of arithmetic genus $25$. The flow can be integrated explicitly in terms of the function $\Upsilon$ such that $(\Upsilon,\Upsilon')$ parametrizes a certain curve of arithmetic genus $3$. 
\end{itemize}
All three vector fields are solenoidal, the first one is also irrotational. Any other $3$-dimensional irreducible superflow is linearly conjugate to one of these. 
\end{thm}
\begin{thm}
\label{thm2}
There exist the following two reducible projective superflows over $\mathbb{R}$ in dimension $3$. 
\begin{itemize}
\item[i)] The superflow $\phi_{\mathbb{P}_{3}}$ with the vector field
\begin{eqnarray*}
(-x^2+2xy+y^2)\frac{\p}{\p x}+(x^2+2xy-y^2)\frac{\p}{\p y}+0\frac{\p}{\p z},
\end{eqnarray*}
whose group of symmetries is the group of order $12$, all symmetries of a $3$-prism (as a group, isomorphic to a dihedral group $\mathbb{D}_{6}$). The orbits are space curves $\{(x-y)(x^2+4xy+y^2)=\mathrm{const.}, z=\mathrm{const.}\}$, and generically these are elliptic curves. The flow can be integrated explicitly in terms of Dixonian elliptic functions.
\item[ii)] The  superflow $\phi_{\mathbb{A}_{4}}$ with the vector field
\begin{eqnarray*}
\frac{x^3-3xy^2}{z}\frac{\p}{\p x}+\frac{y^3-3yx^2}{z}\frac{\p}{\p y}+0\frac{\p}{\p z},
\end{eqnarray*}
whose group of symmetries is the group of order $16$, all symmetries of a $4$-antiprism (as a group, isomorphic to a dihedral group $\mathbb{D}_{8}$). The orbits are curves $\{x^3y-xy^3=\mathrm{const.}, z=\mathrm{const.}\}$, and generically these are curves of arithmetic genus $3$.
\end{itemize}
Both superflows are solenoidal. Any other $3$-dimensional reducible superflow is linearly conjugate to one of these.
\end{thm}

\begin{cor}
\label{cor}
In dimension $3$ over $\mathbb{R}$, there exist a finite number of (irreducible as well as reducible) superflows, and they all are solenoidal.
\end{cor}
In \cite{alkauskas-un,alkauskas-ab,alkauskas-comm,alkauskas-super1} we presented a list of open problems related to the projective translation equation. Concerning the Corrolary, we have seen in \cite{alkauskas-super1} that dihedral superflows in dimension $2$ for groups $\mathbb{D}_{2d+1}$, $d\geq 2$ are not solenoidal. In this aspect, we wonder is this dimension an exception. That is, as a refinement of a question asked in \cite{alkauskas-super1}, we may pose
\begin{prob12}Does the statement of the Corollary hold in any dimension $n\geq 3$? 
\end{prob12}
Define now the numbers 
\begin{eqnarray*}
\psi_{n}=\inf\limits_{\Gamma}\frac{|\Gamma|}{n!},\quad
\Psi_{n}=\sup\limits_{\Gamma}\frac{|\Gamma|}{(n+1)!}.
\end{eqnarray*}
Infimum and supremum are taken over all finite subgroups of $O(n)$ for which there exist a superflow. For example, Theorems \ref{thm1} and \ref{thm2} imply that
\begin{eqnarray*}
\psi_{3}=2,\quad \Psi_{3}=\frac{5}{2}.
\end{eqnarray*} 
For each $n\geq 2$, in \cite{alkauskas-un,alkauskas-super1} we constructed
a superflow without denominators (that is, given by a collection of $n$ quadratic forms) such that it is a superflow for a standard $n$-dimensional representation of a symmetric group $S_{n+1}$. These superflows are solenoidal. In dimension $3$ this is conjugate to the tetrahedral superflow $\phi_{\widehat{\mathbb{T}}}$. Thus,
\begin{eqnarray*}
\Psi_{n}\geq 1\text{ for }n\geq 3.
\end{eqnarray*}
Theorem \ref{thm2} describes the reducible superflow $\phi_{\mathbb{P}_{3}}$, which is just a direct sum of the dihedral $2$-dimensional superflow $\phi_{\mathbb{D}_{3}}$ and the trivial $1$-dimensional flow $\phi(x)=x$. However, the group of symmetries is that of the $3$-prism, and is isomorphic to $\mathbb{P}_{3}\simeq
\mathbb{D}_{3}\oplus\mathbb{Z}_{2}$. In this direction, we prove the following.
\begin{prop}For each $n\geq 3$, there exist a reducible superflow $\Delta_{n}$ of dimension $n$ whose group of symmetries $\Gamma$ is of order $2n!$. As a group, is isomorphic to $S_{n}\oplus\mathbb{Z}_{2}$. 
\label{prop-ext}
\end{prop} 
Of course, we are interested in symmetry groups up to equivalency as representations, which is a finer notion than just an isomorphism. In the proof the exact represention is given.\\
   
In particular, Proposition \ref{prop-ext} implies
\begin{eqnarray*}
\psi_{n}\leq 2\text{ for }n\geq 3.
\end{eqnarray*}
Dimension $2$ is yet again an exception, since
\begin{eqnarray*}
\psi_{2}=3,\quad \Psi_{2}=\infty.
\end{eqnarray*} 
\begin{prob13}Is it true that $\psi_{n}=2$ for $n\geq 3$?
\end{prob13}

\section{Proofs of main results}
\begin{proof} \emph{(Theorem \ref{thm1}, Part I)}. The first two superflows were shown to exist in \cite{alkauskas-super1}. The third superflow is constructed in Section \ref{sec-ico}. There is one more point to explain. The group $A_{5}$ (isomorphic to the icosahedral group) has two irreducible exact $3$-dimensional representations. Let $\tau:\mathbb{Q}\big{[}\sqrt{5}\,\big{]}\mapsto\mathbb{Q}\big{[}\sqrt{5}\,\big{]}$ be a non-trivial automorphism of a quadratic number field. Then this second embedding of an icosahedral group is generated by matrices $\tau(\alpha)=\alpha$, $\tau(\beta)=\beta$ and $\tau(\gamma)\neq\gamma$ (see Subsection \ref{inv-ve}), and these two representations are non-equivalent \cite{kostrikin}. Indeed, if we consider $\mathbb{I}$ as an alternating group $A_{5}\triangleleft S_{5}$, then it possesses an outer automorphism $\pi\mapsto (12)^{-1}\pi(12)$, $\pi\in A_{5}$, which interchanges two different conjugacy classes with representatives $(12345)$ and $(12354)$, thus interchanging character values $\frac{1+\sqrt{5}}{2}$ and $\frac{1-\sqrt{5}}{2}$. Therefore, there exists another superflow with the vector field 
$\tau(\varpi)\bl\tau(\varrho)\bl\tau(\sigma)$. However, let
\begin{eqnarray*}
\epsilon=\begin{pmatrix}
0 & 1 & 0\\
1 & 0 & 0\\
0 & 0 & 1 
\end{pmatrix}.
\end{eqnarray*}
Then it is easy to see that $\epsilon^{-1}\circ(\varpi\bl\varrho\bl\sigma)\circ\epsilon=\tau(\varpi)\bl\tau(\varrho)\bl\tau(\sigma)$, and thus these two superflows are linearly conjugate. \\

The superflows with groups \textbf{(1)}, \textbf{(2)}, and \textbf{(3)} exist.\\

The groups given as direct product with $\{I,-I\}$ do not produce superflows. Indeed, as we already know, $-I$ cannot be an element for the symmetry group for the superflow, since (\ref{kappa}) shows that then the vector field is trivial. The rest of the groups are reducible $3$-dimensional representations, so this finishes the proof of Theorem \ref{thm1}.
\end{proof}

\begin{proof} \emph{(Theorem \ref{thm2})}. Equally, as noted in \cite{alkauskas-super1}, the group $\mathbb{C}_{n}$ does not produce a superflow, too, since the invariant vector field is never unique. See \textbf{(4)} below for an argument, which is completely analogous in a rotation group case. The group $\mathbb{T}$ produces a superflow, but its group of symmetries extends to $\widehat{\mathbb{T}}$. This phenomenon is explained in (Note 5, \cite{alkauskas-super1}). Now we pass to the rest $4$ cases.\\ 

\indent\textbf{(4)}. Let $d\in\mathbb{N}$. The mixed cyclic group $\mathbb{C}_{2d}\mathbb{C}_{d}$ of order $2d$ is generated by a matrix
\begin{eqnarray*}
\alpha=\begin{pmatrix}
-\cos\big{(}\frac{\pi}{d}\big{)} & \sin\big{(}\frac{\pi}{d}\big{)} & 0\\
-\sin\big{(}\frac{\pi}{d}\big{)} & -\cos\big{(}\frac{\pi}{d}\big{)} & 0\\
0 & 0 & -1 
\end{pmatrix}.
\end{eqnarray*}
Let
\begin{eqnarray}
\tau=\begin{pmatrix}
\frac{1}{\sqrt{2}} & \frac{i}{\sqrt{2}} & 0\\
\frac{i}{\sqrt{2}} & \frac{1}{\sqrt{2}} & 0\\
0 & 0 & 1 
\end{pmatrix}.\label{tau}
\end{eqnarray} 
By a direct calculation,
\begin{eqnarray*}
\hat{\alpha}=\tau\cdot\alpha\cdot\tau^{-1}=\begin{pmatrix}
-\zeta & 0 & 0\\
0 & -\zeta^{-1} & 0\\
0 & 0 & -1 
\end{pmatrix},\quad\zeta=e^{\frac{\pi i}{d}}.
\end{eqnarray*}
Suppose, a vector field (over $\mathbb{C}$) $\varpi\bl\varrho\bl\sigma$ is invariant under conjugation with $\hat{\alpha}$. Then so is 
\begin{eqnarray}
\varpi(x,y,z)\bl a\varpi(y,x,z)\bl b\sigma(x,y,z),\quad a,b,c\in\mathbb{R}.
\label{3-pam}
\end{eqnarray}
Indeed, by the assumption, 
\begin{eqnarray}
-\zeta^{-1}\varpi\big{(}-\zeta x,-\zeta^{-1}y,-z\big{)}=\varpi(x,y,z).
\label{sup}
\end{eqnarray} 
Therefore,
\begin{eqnarray*}
-\zeta\varpi\big{(}-\zeta^{-1}y,-\zeta x,-z\big{)}
=-\zeta^{-2d+1}
\varpi\big{(}-\zeta^{2d-1}y,-\zeta^{-2d+1}x,-z\big{)}\mathop{=}^{(\ref{sup})}
\varpi(y,x,z).
\end{eqnarray*}
Since at least one of $\varpi$ and $\sigma$ is not identically zero, the vector field (\ref{3-pam}) is at least a $1$ parameter family, and thus $\alpha$ does not produce a superflow.\\

\indent\textbf{(5)}. Let $d\in\mathbb{N}$. The dihedral group is generated by two matrices
\begin{eqnarray}
\begin{pmatrix}
\cos\big{(}\frac{2\pi}{d}\big{)} & -\sin\big{(}\frac{2\pi}{d}\big{)} & 0\\
\sin\big{(}\frac{2\pi}{d}\big{)} & \cos\big{(}\frac{2\pi}{d}\big{)} & 0\\
0 & 0 & 1 
\end{pmatrix},\quad
\epsilon=\begin{pmatrix}
0 & 1 & 0\\
1 & 0 &0\\
0 & 0 & 1 
\end{pmatrix}.
\label{epsi}
\end{eqnarray}
 This does not produce a superflow. Indeed, the third coordinate split, and if $\varpi\bl\varrho\bl\sigma$ is an invariant vector field, so is a one parameter family $\varpi\bl\varrho\bl az^2$, $a\in\mathbb{R}$. This contradicts the definition of a superflow.\\

\indent\textbf{(6)}. This time the group $\Gamma$ isomorphic to a dihedral group $\mathbb{D}_{d}$ realized as a group of rotations, and it is generated by two matrices
\begin{eqnarray}
\alpha=\begin{pmatrix}
\cos\big{(}\frac{2\pi}{d}\big{)} & -\sin\big{(}\frac{2\pi}{d}\big{)} & 0\\
\sin\big{(}\frac{2\pi}{d}\big{)} & \cos\big{(}\frac{2\pi}{d}\big{)} & 0\\
0 & 0 & 1 
\end{pmatrix},\quad
\beta=\begin{pmatrix}
0 & 1 & 0\\
1 & 0 &0\\
0 & 0 & -1 
\end{pmatrix}
\label{dih}
\end{eqnarray}
This case needs more analysis. Let $\tau$ be given by (\ref{tau}). Let $\zeta=e^{\frac{2\pi i}{d}}$. As calculations show, 
\begin{eqnarray*}
\hat{\alpha}=\tau\alpha\tau^{-1}=\begin{pmatrix}
\zeta & 0 & 0\\
0 & \zeta^{-1} & 0\\
0 & 0 & 1 
\end{pmatrix},\quad
\hat{\beta}=\tau\beta\tau^{-1}=\begin{pmatrix}
0 & 1 & 0\\
1 & 0 & 0\\
0 & 0 & -1 
\end{pmatrix}.
\end{eqnarray*}
So, we need to find a vector field over $\mathbb{C}$ which is invariant under conjugation with these two matrices.\\

We can finish this case more quickly due to the following argument. As we know, the denominator of the superflow is a relative invariant of the group \cite{alkauskas-super1}. The group generated by $\hat{\alpha}$ and $\hat{\beta}$ has two invariants of degree $2$; namely, $xy$ and $z^2$. Thus, the denominator for the superflow for this group cannot be of degree higher than one, otherwise it is not unique - we can replace the appearance of $xy$ by $z^2$, and vice versa. We are left to consider the case when it is of degree $0$ or $1$ (then the denominator is necessarily $z$, a relative invariant). \\

When denominator is of degree $0$, with some experimentation we get that $d=3$, and there exists the unique (up to homothety) vector field $y^2\bl x^2\bl 0$ with the symmetry group of order $6$. As we will see soon, its group of symmetries extends to the group of order $12$: to matrices $\alpha$ and $\beta$ we may add the matrix $\epsilon$, as given by (\ref{epsi}).\\

When the denominator is $z$, we find that $d=4$ and the invariant vector field is given by
\begin{eqnarray*}
\frac{y^3}{z}\bl-\frac{x^3}{z}\bl 0.
\end{eqnarray*}
Again, to its group of symmetries we may add the matrix $\hat{\gamma}$ as given by (\ref{neg}) (where $\xi=e^{\frac{\pi i }{4}}$, $\xi^{4}=-1$), and so the group of symmetries, the one of order $8$, extends to a group of order $16$. In both cases, we may pass to the last case left in the Table \ref{table1}.\\
 
\indent\textbf{(7).} We will expand on this case more extensively. Let $\ell\in\mathbb{N}$. Consider the mixed dihedral group $\widehat{\Gamma}=\mathbb{D}_{2\ell}\mathbb{D}_{\ell}$. If $\ell$ is odd, then this is the full symmetry group of an $\ell$-prism - the Cartesian product of a regular $\ell$-gon in the horizontal plane with an interval $[-1,1]$ in the vertical plane. If $\ell$ is even, $\widehat{\Gamma}$ is the full symmetry group of the so called \emph{antiprism}. The latter is made of two regular $\ell$-gons in planes $z=-1$ and $z=1$, but one turned by an angle of $\frac{\pi}{\ell}$ with respect to the other, and then connected with isosceles triangles. See Figure \ref{anti} for an image of a square antiprism. \\

\begin{figure}
\epsfig{file=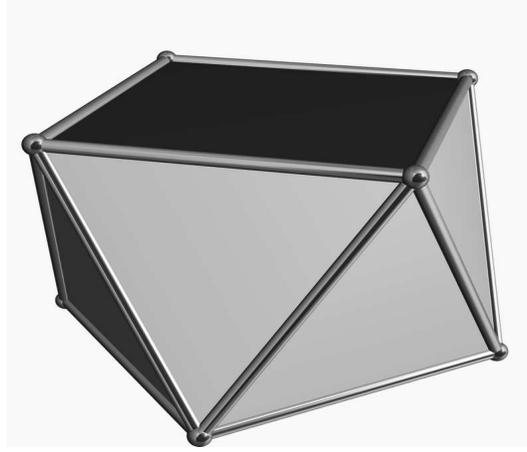,width=200pt,height=170pt,angle=0}
\caption{The square antiprism}
\label{anti}
\end{figure}

Suppose now $\ell$ is even, $\ell=2d$, $d\in\mathbb{N}$. Let $\xi=e^{\frac{\pi i}{\ell}}$. The group $\widehat{\Gamma}$ is generated by $\gamma$ and $\epsilon$, where matrices $\gamma$ and $\epsilon$ are given by
\begin{eqnarray}
\gamma=\begin{pmatrix}
\cos\big{(}\frac{\pi}{\ell}\big{)} &-\sin\big{(}\frac{\pi}{\ell}\big{)} & 0\\
\sin\big{(}\frac{\pi}{\ell}\big{)} &\cos\big{(}\frac{\pi}{\ell}\big{)} & 0\\
0 & 0 & -1 
\end{pmatrix},\quad
\epsilon=\begin{pmatrix}
0 & 1 & 0\\
1 & 0 & 0\\
0 & 0 & 1 
\end{pmatrix}.
\label{epsi2}
\end{eqnarray}
If $\tau$ is again given by (\ref{tau}), then
\begin{eqnarray}
\hat{\gamma}=\tau\gamma\tau^{-1}=\begin{pmatrix}
\xi & 0 & 0\\
0 & \xi^{-1} & 0\\
0 & 0 & -1 
\end{pmatrix},\quad
\hat{\epsilon}=\tau\epsilon\tau^{-1}=\begin{pmatrix}
0 & 1 & 0\\
1 & 0 & 0\\
0 & 0 & 1 
\end{pmatrix}.
\label{neg}
\end{eqnarray}
 Thus, $\xi^{2d}=\xi^{\ell}=-1$. Suppose $d\geq 2$. We can argue similarly as in the case \textbf{(6)} to claim that the degree of the denominator should be equal to $0$ or $1$. Here we present a slightly different apporach. We verify directly that there exists no invariant vector field with degree of denominator $<2d-3$, and for degree exactly $2d-3$ the family of invariant vector fields is given by
\begin{eqnarray}
\mathbf{C}_{a}(\m{x})=\varpi\bl\varrho\bl\sigma=\frac{y^{2d-1}}{(xy)^az^{2d-3-2a}}\bl\frac{x^{2d-1}}{(xy)^az^{2d-3-2a}}\bl 0,\quad 0\leq a\leq d-2.\label{ca}
\end{eqnarray}
In general, the following $(d-1)$-parameter family is an invariant vector field:
\begin{eqnarray*}
\mathbf{C}(\m{x})=\frac{y^{2d-1}}{\sum\limits_{a=0}^{d-2}c_{a}(xy)^az^{2d-3-2a}}\bl\frac{x^{2d-1}}{\sum\limits_{a=0}^{d-2}c_{a}(xy)^az^{2d-3-2a}}\bl 0,\quad 0\leq a\leq d-2.
\end{eqnarray*}
Up to homothety, this is a $(d-2)$-parameter family. Note that (\ref{ca}) for $d=1$ does not work since we have a $1$-parameter family of invariant vector fields (in this case, $\xi=i$)
\begin{eqnarray*}
yz\bl xz\bl c(x^{2}+y^2).
\end{eqnarray*}

Invariance of $\mathbf{C}_{a}$ under $\hat{\epsilon}$ is immediate, and
\begin{eqnarray*}
\hat{\gamma}^{-1}\circ\mathbf{C}_{a}\circ\hat{\gamma}(\m{x})
&=&\hat{\gamma}^{-1}\circ(\varpi\bl\varrho\bl\sigma)\circ\hat{\gamma}(x,y,z)\\
&=&\hat{\gamma}^{-1}\circ(\varpi\bl\varrho\bl\sigma)\circ(\xi x,\xi^{-1}y,-z)\\
&=&\hat{\gamma}^{-1}\circ\Big{(}\frac{\xi^{-2d+1}y^{2d-1}}{(xy)^a(-z)^{2d-3-2a}}
\bl\frac{\xi^{2d-1}x^{2d-1}}{(xy)^a(-z)^{2d-3-2a}}\bl 0\Big{)}\\
&=&-\xi^{-1}\frac{\xi^{-2d+1}y^{2d-1}}{(xy)^az^{2d-3-2a}}
\bl-\xi\frac{\xi^{2d-1}x^{2d-1}}{(xy)^az^{2d-3-2a}}\bl 0\\
&=&-\xi^{-2d}\frac{y^{2d-1}}{(xy)^az^{2d-3-2a}}\bl-\xi^{2d}\frac{x^{2d-1}}{(xy)^az^{2d-3-2a}}\bl 0\\
&=&\mathbf{C}_{a}(\m{x}).
\end{eqnarray*}
It is obvious now that the vector field $\mathbf{C}(\m{x})$ is also an invariant vector field. This up to homothety is defined uniquely only if $d-2=0$. Then it is equal to $\frac{y^3}{z}\bl\frac{x^3}{z}\bl 0$. Note that this vector field is conjugate to the vector field $\frac{y^3}{z}\bl-\frac{x^3}{z}\bl 0$ we obtained in part \textbf{(6)}. This is seen with a help of conjugation with a diagonal matrix $\mathrm{diag}(\xi^2,\xi,\xi)$, $\xi^{4}=-1$. Thus, if $d=2$, we  get a reducible superflow! Let us return to the real setting, and calculate the vector field $\tau^{-1}\circ\mathbf{C}\circ\tau$. We thus obtain the vector field
\begin{eqnarray}
\frac{x^3-3xy^2}{z}\bl \frac{y^3-3yx^2}{z}\bl 0,
\label{anti-4}
\end{eqnarray}
with orbits $\{\mathscr{W}(x,y)=x^3y-xy^3=\mathrm{const.}, z=\mathrm{const}\}$. In this manner we get the superflow $\phi_{\mathbb{A}_{4}}$, which has a group of a $4$-antiprism of order $16$ as its symmetry group: it is generated by two matrices
\begin{eqnarray*}
\begin{pmatrix}
\frac{\sqrt{2}}{2} &-\frac{\sqrt{2}}{2} & 0\\
\frac{\sqrt{2}}{2} &\frac{\sqrt{2}}{2} & 0\\
0 & 0 & -1 
\end{pmatrix}\text{ and }
\begin{pmatrix}
0 & 1 & 0\\
1 & 0 & 0\\
0 & 0 & 1 
\end{pmatrix}.
\end{eqnarray*}

 Figure \ref{quadratic} explains the beautiful visualisation of an antiprism that we get this way. This vector field and the superflow itself is explored in Section \ref{sec-anti}.\\
 
\begin{figure}
\epsfig{file=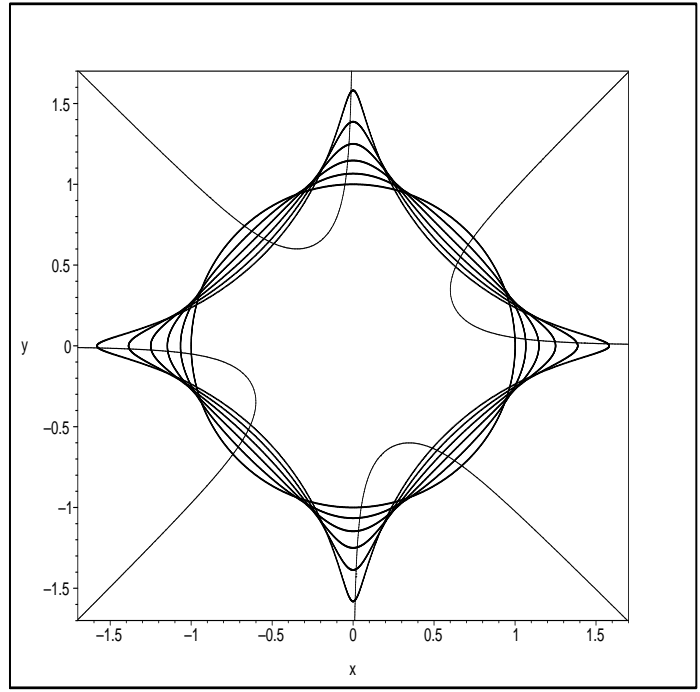,width=260pt,height=260pt,angle=-90}
\caption{Deformation of the unit circle $\mathbf{S}=\{x^{2}+y^{2}=1\}$ in the plane $z=1$ under the antiprismal superflow $\phi_{\mathbb{A}_{4}}$. The seven curves $\phi_{\mathbb{A}_{4}}^{0.06j}(\mathbf{S})$, $j=0,\ldots 6$, are shown. Since the vector field is solenoidal, the area inside each curve is equal to $\pi$. The selected orbit $x^3y-xy^3=0.05$ is also shown. Note that in \cite{alkauskas-super1} we presented pictures of a $5$-point star and $7$-point star, both arising as visualizations of dihedral irreducible superflows. The appearance of a $4$-point star is new, and could not be realized in a $2$-dimensional setting. Further, in the plane $z=-1$ the picture is analogous, though different - we get the same vector field, but with direction of the flow inverted, as is clear from (\ref{anti-4}). Therefore, in a plane $z=-1$ we get the same $4$-point star, but turned by the angle $\frac{\pi}{4}$; see Figure \ref{anti}. Now the appearance of a $4$-antiprism symmetry group as symmetries of the flow $\phi_{\mathbb{A}_{4}}$ becomes evident. When $z$ varies from $1$ to $0$, the deformation becomes larger. Thus, the cylinder $x^2+y^2$ under this flow transforms like square antiprism, and the picture above might also be interpreted as a deformation of the  infinite cylinder $\{x^2+y^2=1,1\leq z<\infty\}$ at a fixed time moment, seen from below.}
\label{quadratic}
\end{figure}

Now suppose that $\ell=2d+1$ is odd. The matrices $\gamma$ and $\epsilon$ still produce the symmetry group of a $\ell$-antiprism, but this time the group is not of a mixed type, but a direct product as given by Table \ref{table1}
in the middle column, since $\gamma^{\ell}=-I$. As said, we need to consider the symmetry group of a $\ell$-prism, and therefore the group $\widehat{\Gamma}=\mathbb{D}_{2\ell}\mathbb{D}_{\ell}$ is generated by
 \begin{eqnarray*}
\gamma=\begin{pmatrix}
\cos\big{(}\frac{2\pi}{\ell}\big{)} &-\sin\big{(}\frac{2\pi}{\ell}\big{)} & 0\\
\sin\big{(}\frac{2\pi}{\ell}\big{)} &\cos\big{(}\frac{2\pi}{\ell}\big{)} & 0\\
0 & 0 & -1 
\end{pmatrix},\quad
\epsilon=\begin{pmatrix}
0 & 1 & 0\\
1 & 0 & 0\\
0 & 0 & 1 
\end{pmatrix}.
\end{eqnarray*}
Again, conjugating with $\tau$ we obtain matrices $\hat{\gamma}$ and $\hat{\epsilon}$ as given by (\ref{neg}), but this time $\xi=e^{\frac{2\pi i}{\ell}}$. We verify that there exist no invariant vector field with denominator of degree lower than $2d-2$, and for degree exactly $2d-2$ the family of invariant vector fields is given by
\begin{eqnarray*}
\mathbf{C}(\m{x})=\varpi\bl\varrho\bl\sigma=\frac{y^{2d}}{\sum\limits_{a=0}^{d-1}c_{a}(xy)^az^{2d-2-2a}}\bl\frac{x^{2d}}{\sum\limits_{a=0}^{d-1}c_{a}(xy)^az^{2d-2-2a}}\bl 0,\quad 0\leq a\leq d-1.
\end{eqnarray*} 
The vector field (up to homothety, of course) is unique only if $d=1$. And indeed, for $\xi=e^{\frac{2\pi i}{3}}$ the unique invariant vector field, which has no denominator, is given by $y^2\bl x^2\bl 0$. Returning to the real case, we obtain the vector field
\begin{eqnarray*}
\varpi\bl\varrho\bl\sigma=-x^2+2xy+y^2\bl x^2+2xy-y^2\bl 0.
\end{eqnarray*}
In this manner we get the superflow $\phi_{\mathbb{P}_{3}}$, which has a group of a $3$-prism of order $12$ as its symmetry group: it is generated by two matrices
\begin{eqnarray*}
\begin{pmatrix}
-\frac{1}{2} &-\frac{\sqrt{3}}{2} & 0\\
\frac{\sqrt{3}}{2} &-\frac{1}{2} & 0\\
0 & 0 & -1 
\end{pmatrix}\text{ and }
\begin{pmatrix}
0 & 1 & 0\\
1 & 0 & 0\\
0 & 0 & 1 
\end{pmatrix}.
\end{eqnarray*}
In fact, we do not need to integrate this vector field - if we perform a linear conjugation with a matrix 
\begin{eqnarray*}
\delta=\begin{pmatrix}
1-\sqrt{3} &  1+\sqrt{3}& 0\\
1+\sqrt{3}&1-\sqrt{3} & 0\\
0 & 0 & 1 
\end{pmatrix},
\end{eqnarray*}
then
\begin{eqnarray*}
\delta^{-1}\circ(\varpi\bl\varrho\bl\sigma)\circ\delta=-4(x^2-2xy)\bl-4(y^2-2xy)\bl 0.
\end{eqnarray*}
The vector field $x^2-2xy\bl y^2-2xy$ and the superflow itself was extensively studied in \cite{alkauskas-un}. We only choose to use the vector field $\varpi\bl\varrho\bl\sigma$ in the formulation of case i) of Theorem \ref{thm2} only because its symmetry group is inside $O(3)$ (but defined over $\mathbb{Q}[\sqrt{3}\,]$), while that of $x^2-2xy\bl y^2-2xy\bl 0$ is inside $\mathrm{GL}(2,\mathbb{R})$ but not $O(3)$, though defined over $\mathbb{Q}$.  
\end{proof}

\begin{proof}\emph{ (Proposition \ref{prop-ext})}. Here we reproduce the result proved in  \cite{alkauskas-un}, also mentioned in \cite{alkauskas-super1}.\\

Let $n\in\mathbb{N}$, $n\geq 2$, and consider the following vector field
\begin{eqnarray}
\mathbf{Q}(\m{x})=Q_{1}(\mathbf{x})\bl Q_{2}(\mathbf{x})\bl \ldots\bl Q_{n}(\mathbf{x}).
\label{Q-vec}
\end{eqnarray}
Here
\begin{eqnarray*}
Q_{1}(\m{x})=x_{1}^{2}-\frac{2}{n-1}\cdot x_{1}
\sum\limits_{i=2}^{n}x_{i},\quad \m{x}=x_{1}\bl x_{2}\bl\ldots\bl x_{n},
\end{eqnarray*}
and $Q_{i}$ are obtained from the above component by interchanging the r\^{o}les of $x_{1}$ and $x_{i}$. Now, consider the standard permutation representation of the symmetric group $S_{n}$. It is clear that the vector field $\mathbf{Q}$ is invariant under conjugation with all elements of $S_{n}$. Moreover, consider the following $n\times n$ matrix of order $2$:
\begin{eqnarray*}
\kappa=\begin{pmatrix}
-1 &   &   &   &  \\
-1 & 1 &   &   &  \\
-1 &   & 1 &   &  \\
\,\,\,\vdots & &   & \ddots &  \\
-1 &   &  &  & 1
\end{pmatrix}.
\end{eqnarray*} 
Then (\ref{kappa}) is satisfied, and this is easily verified by hand. The group $S_{n}$ and the involution $\kappa$ generate the group $S_{n+1}$. This is a vector field of a $n$-dimensional irreducible superflow. Let $\sigma:S_{n+1}\mapsto\mathrm{GL}(n,\mathbb{Q})\subset\mathrm{GL}(n,\mathbb{R})$ be this exact representation of a group $S_{n+1}$.\\

Now (here the new material starts), consider the group $S_{n+1}\oplus\mathbb{Z}_{2}$, where $\mathbb{Z}_{2}=\{1,-1\}$ as a multiplicative group. Let $(\tau,\epsilon)\in S_{n+1}\oplus\mathbb{Z}_{2}$. Consider the $(n+1)$-dimensional representation
\begin{eqnarray*}
\widehat{\sigma}:S_{n+1}\oplus\mathbb{Z}_{2}\mapsto\mathrm{GL}(n+1,\mathbb{R}),\quad \widehat{\sigma}\big{(}(\tau,\epsilon)\big{)}=\sigma(\tau)\oplus\tilde{\epsilon};
\end{eqnarray*}
here $\tilde{\epsilon}$ is a $1\times 1$ scalar matrix with an entry $\epsilon$. Then the representaion is exact, and let $\Gamma=\widehat{\sigma}(S_{n+1}\oplus\mathbb{Z}_{2})$. We will show that there exist a superflow for $\Gamma$, and its vector field is given by
\begin{eqnarray*}
\widehat{\mathbf{Q}}(\hat{\m{x}})=Q_{1}(\mathbf{x})\bl Q_{2}(\mathbf{x})\bl \ldots\bl Q_{n}(\mathbf{x})\bl 0,\quad 
\hat{\m{x}}=x_{1}\bl x_{2}\bl\ldots\bl x_{n+1}.
\end{eqnarray*}
Of course, this is just a generalization of the reducible superflow in Theorem \ref{thm2}, case i). Thus, the vector field $\widehat{\mathbf{Q}}$ is independend of the coordinate $x_{n+1}$.\\

Suppose now that a vector field $\mathbf{R}$, a collection of $(n+1)$-quadratic forms, is invariant under conjugation with all elements from $\Gamma$. Since the matrix $\chi=\mathrm{diag}(\underbrace{1,\ldots,1}_{n},-1)$ belongs to $\Gamma$, we readily obtain that the first $n$ components of $\mathbf{R}$ are independent of $x_{n+1}$, except possibly for a summand proportional to $x_{n+1}^{2}$. Since $\widehat{\sigma}(S_{n+1}\oplus\{1\})$ is a subgroup of $\Gamma$, we get that the first $n$ components are given by a scalar multiple of $\mathbf{Q}$ as given by (\ref{Q-vec}), plus this extra summand. Again, using (\ref{kappa}) for a matrix $\chi$, we see that the invariant vector field for the group $\Gamma$ is necessarily of the form
\begin{eqnarray*}
Q_{1}+a_{1}x_{n+1}^{2}\bl Q_{2}+a_{2}x_{n+1}^{2}\bl\ldots\bl Q_{n}+a_{n}x_{n+1}^{2}\bl x_{n+1}(\sum\limits_{j=1}^{n}b_{j}x_{j}).
\end{eqnarray*}
Now, $S_{n}\times\{1\}$ is a subgroup of $\Gamma$, since we constructed the representation of $S_{n+1}$ from a permutation representation of $S_{n}$ adding one extra order two matrix $\kappa$. Checking (\ref{kappa}) for $\gamma\in S_{n}\times\{1\}$ we get that the vector field is actually
\begin{eqnarray*}
Q_{1}+ax_{n+1}^{2}\bl Q_{2}+ax_{n+1}^{2}\bl\ldots\bl Q_{n}+ax_{n+1}^{2}\bl bx_{n+1}(\sum\limits_{j=1}^{n}x_{j}),
\quad a,b\in\mathbb{R}. 
\end{eqnarray*}
We can finish now by checking the condition (\ref{kappa}) for the matrix $\gamma=\kappa\oplus\{1\}$. We get that $a=b=0$. This proves Proposition \ref{prop-ext} (with $n$ in place of $n+1$, if compared to the proof above).
\end{proof}
\section{Icosahedral superflow}
\label{sec-ico}
\subsection{Invariant vector field}
\label{inv-ve}
Let us define
\begin{eqnarray*}
\alpha\mapsto\begin{pmatrix}
-1 & 0 & 0\\
0 & -1 & 0\\
0 & 0 & 1 
\end{pmatrix},\quad
\beta\mapsto\begin{pmatrix}
0 & 0 & 1\\
1 & 0 &0\\
0 & 1 & 0 
\end{pmatrix},\quad
\gamma\mapsto\begin{pmatrix}
\frac{1}{2} & -\frac{\phi}{2} & \frac{1}{2\phi}\\
\frac{\phi}{2} & \frac{1}{2\phi} & -\frac{1}{2}\\
\frac{1}{2\phi} & \frac{1}{2} & \frac{\phi}{2} 
\end{pmatrix}.
\end{eqnarray*}
These are matrices of orders $\alpha^2=\beta^3=\gamma^{5}=I$, and together they generate the icosahedral group $I_{h}$ of order $60$. Let us call it $\mathbb{I}\subset SO(3)$.\\


As is implied from \cite{goller}, the group $\mathbb{I}$ has a unique, up to scalar multiple, invariant of degree $4$, and it is $D=(x^2+y^2+z^2)^2$. The general form of a $2$-homogenic vector field with a denominator $D$, which is invariant under conjugation with matrices $\alpha$ and $\beta$, is given by $\varpi\bl\varrho\bl\sigma$, where
\begin{eqnarray*}
\varpi=\frac{Ayz^5+By^5z+Cy^3z^3+Dx^2yz^3+Ex^2y^3z+Fx^4yz}{(x^2+y^2+z^2)^2},\\
\varrho=\frac{Azx^5+Bz^5x+Cz^3x^3+Dy^2zx^3+Ey^2z^3x+Fy^4zx}{(x^2+y^2+z^2)^2},\\
\sigma=\frac{Axy^5+Bx^5y+Cx^3y^3+Dz^2xy^3+Ez^2x^3y+Fz^4xy}{(x^2+y^2+z^2)^2}.
\end{eqnarray*}
If we can force $\varpi\bl\varrho\bl\sigma$ to be invariant under conjugation with $\gamma$, we are done! So, let us use the unspecified coefficients $A$ through $F$, and let us calculate $\widehat{\varpi}\bl\widehat{\varrho}\bl\widehat{\sigma}=\gamma^{-1}\circ(\varpi\bl\varrho\bl\sigma)\circ\gamma(x,y,z)$ with the help of MAPLE. The coefficient at $x^6$ in the numerator of $\widehat{\varpi}$ is then equal to a scalar multiple of $C+D+E+2A+2B+2F$. Since the corresponding coefficient of $\varpi$ is equal to $0$, this gives us one linear relation among six coefficients. Continuing in the same manner, that is, equating the coefficients of $\widehat{\varpi}$ to those of $\varpi$, we arrive at $5$ different linear equations of full rank. Thus we find the unique (up to conjugation by a homothety) vector field which is invariant under conjugation with $\gamma$ as well, and this vector field is given by (\ref{vec-field}). 

\subsection{Orbits}
We can check directly that two independent first integrals of the differential system (\ref{sys-in}), i.e. those that satisfy (\ref{first-int}), are given by
\begin{eqnarray*}
\mathscr{W}=x^2+y^2+z^2,\text{ and }\mathscr{V}=(\phi^2x^2-y^2)(\phi^2 y^2-z^2)(\phi^2z^2-x^2).
\end{eqnarray*}
The first one is what was expected. On the other hand, the function $\mathscr{V}$ is an invariant of the group $\mathbb{I}$, and the first integral should be an invariant; direct calculations confirm the above claim. If we consider the group $\widehat{\mathbb{I}}$ of order $120$, then its ring of invariants is polynomial, it is generated by $\mathscr{W}$, $\mathscr{V}$, and the third invariant \cite{goller}
\begin{eqnarray*}
(x+y+z)(x+y-z)(x-y+z)(x-y-z)(\phi^4x^2-y^2)(\phi^4y^2-z^2)(\phi^4z^2-x^2).
\end{eqnarray*}
This is also an invariant of $\mathbb{I}$. However, contrary to $\widehat{\mathbb{I}}$, $\mathbb{I}$ is not generated by pseudo-reflections, so there exists a fourth invariant (Theorem of Chevalley-Shephard-Todd, see \cite{alkauskas-super1,chevalley}).\\

Note that
\begin{eqnarray*}
\varpi_{x}+\varrho_{y}+\sigma_{z}=0,
\end{eqnarray*}
so the vector field is solenoidal. Let this superflow be given by 
\begin{eqnarray*}
\phi(\m{x})=V(x,y,z)\bl V(y,z,x)\bl V(z,x,y).
\end{eqnarray*}

\subsection{Intersection}Consider six planes given by $(\phi^2x^2-y^2)(\phi^2 y^2-z^2)(\phi^2z^2-x^2)=0$. They divide the unit sphere into $32$ regions, and the space into $32$ open solid angles, as follows (``triangles" and ``pentagons" mean regions on the unit sphere bounded by arcs, parts of great circles):
\begin{itemize}
\item[A)]$8$ triangles with centres being  $\frac{1}{\sqrt{3}}(\pm 1,\pm 1,\pm 1)$. Points inside the solid angles are characterized by the fact that all three factors defining $\mathscr{V}$ are positive. The three vertices of the triangle with the center $\frac{1}{\sqrt{3}}(1,1,1)$ are equal to $(\frac{1}{2\phi},\frac{1}{2},\frac{\phi}{2})$, $(\frac{\phi}{2},\frac{1}{2\phi},\frac{1}{2})$, 
$(\frac{1}{2},\frac{\phi}{2},\frac{1}{2\phi})$.

\item[A')]$12$ further triangles; each octant contains three halves.  One of these triangles has a center $\frac{1}{\sqrt{3}}(0,\frac{1}{\phi},\phi)$ and vertices $(\frac{1}{2\phi},\frac{1}{2},\frac{\phi}{2})$, $(0,0,1)$, $(-\frac{1}{2\phi},\frac{1}{2},\frac{\phi}{2})$. Points inside these solid angles are characterized by the fact that two factors defining $\mathscr{V}$ are negative, one is positive. Note, however, that $\gamma(1,1,1)^{T}=(0,\frac{1}{\phi},\phi)$. So, any of these triangles can be identified to any triangle in the list A) with the help of transformation in the group $\mathbb{I}$.  
\item[B)] $12$ pentagons; one of them has a center $\frac{1}{\sqrt{\phi^2+1}}(\phi,1,0)$, and vertices $(1,0,0)$, $(\frac{\phi}{2},\frac{1}{2\phi},\frac{1}{2})$, $(\frac{1}{2},\frac{\phi}{2},\frac{1}{2\phi})$, $(\frac{1}{2},\frac{\phi}{2},-\frac{1}{2\phi})$, $(\frac{\phi}{2},\frac{1}{2\phi},-\frac{1}{2})$.
\end{itemize}

\begin{figure}
\epsfig{file=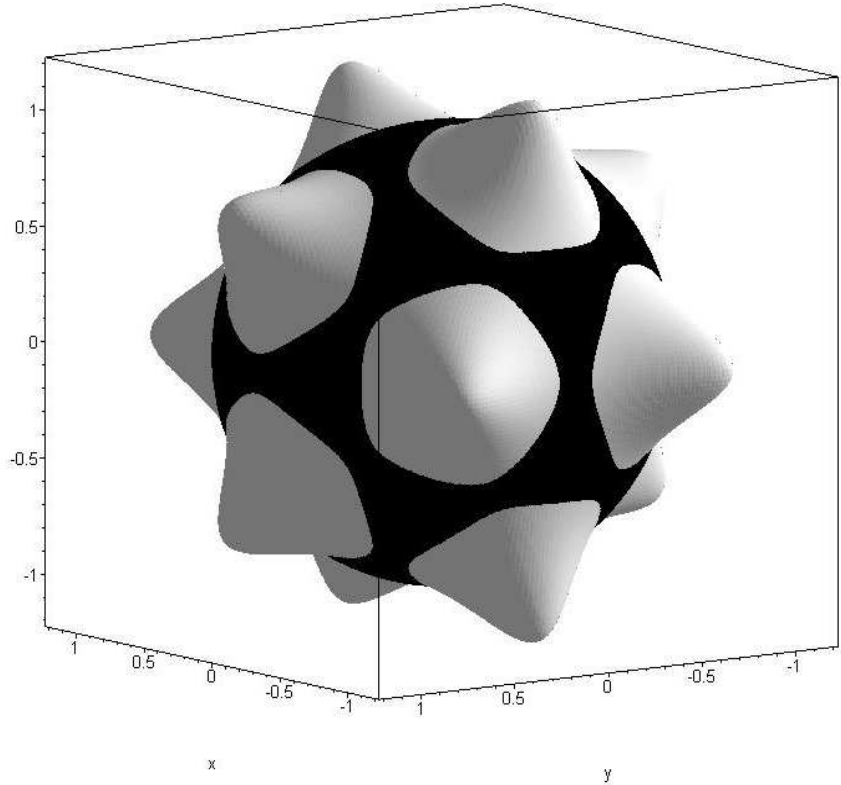,width=380pt,height=350pt,angle=0}
\caption{The intersection of the unit sphere $x^2+y^2+z^2=1$ (black) with the surface  $(\phi^2x^2-y^2)(\phi^2 y^2-z^2)(\phi^2z^2-x^2)+(x^2+y^2+z^2)^3=\frac{19}{20}$ (gray), $\xi=-\frac{1}{20}$.}
\label{figure-neg}
\end{figure}

\begin{figure}
\epsfig{file=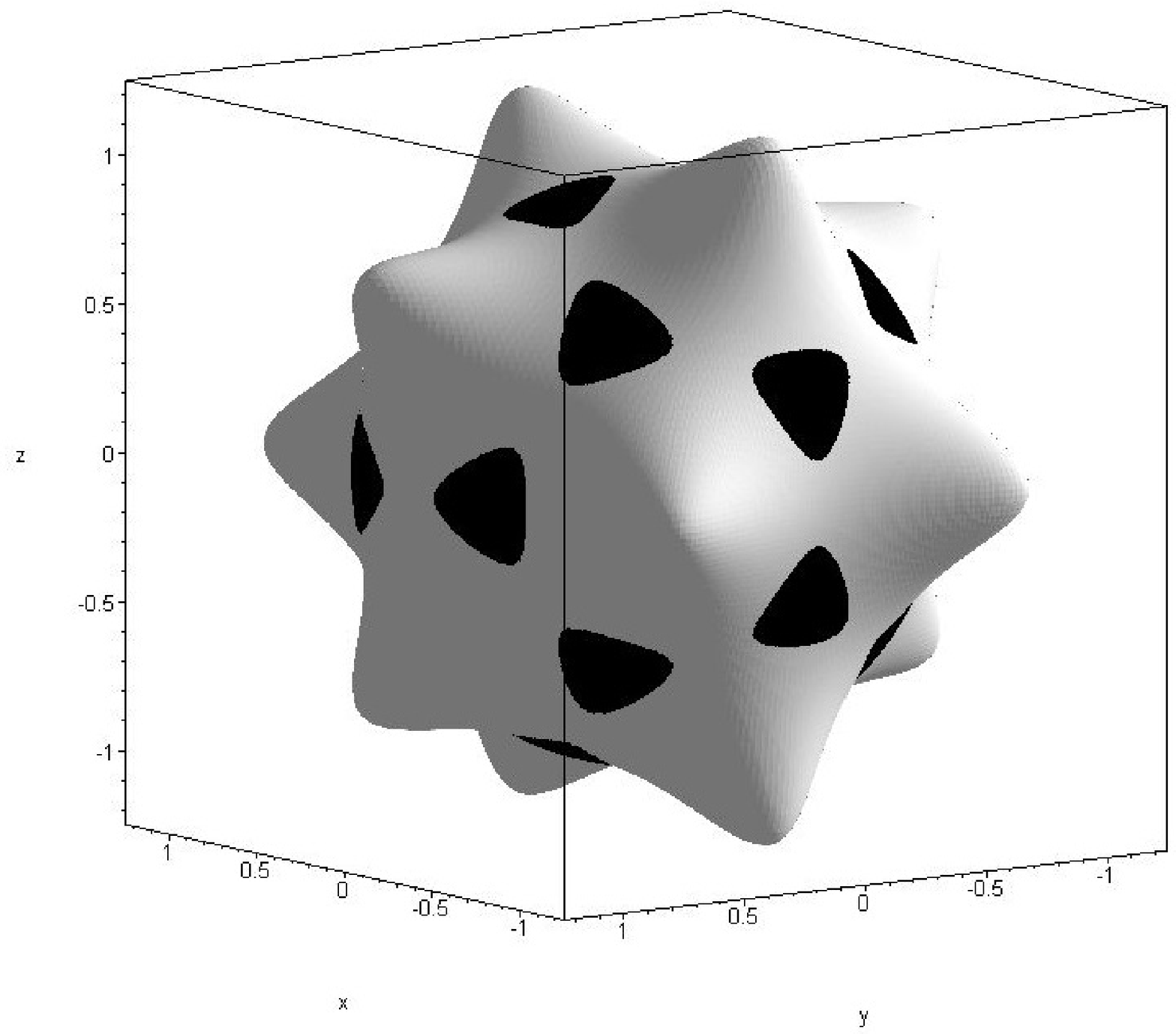,width=380pt,height=350pt,angle=0}
\caption{The intersection of the unit sphere $x^2+y^2+z^2=1$ (black) with the surface  $(\phi^2x^2-y^2)(\phi^2 y^2-z^2)(\phi^2z^2-x^2)+(x^2+y^2+z^2)^3=\frac{21}{20}$ (gray), $\xi=\frac{1}{20}$.}
\label{figure-pos}
\end{figure}

For $\xi\in\mathbb{R}$, let us define the surface $\mathbf{T}_{\xi}=\{\m{x}\in\mathbb{R}^{3}:\mathscr{V}(\m{x})=\xi\}$. 
Thus, we can formulate the following result. 
\begin{prop}Let $\xi\in\mathbb{R}$, and consider the intersection of the unit sphere $\mathbf{S}$ with an algebraic surface $\mathbf{T}_{\xi}$. Then this intersection:
\begin{itemize}
\item[i)] For $\xi<-\frac{\phi^3}{5}=-\frac{2+\sqrt{5}}{5}$, is empty.
\item[ii)]For $\xi=-\frac{2+\sqrt{5}}{5}$, is equal to $12$ isolated points.
\item[iii)]For $-\frac{2+\sqrt{5}}{5}<\xi<0$, consists of $12$ disjoint curves, each homeomorphic to a circle. See (Figure \ref{figure-neg}). 
\item[iv)]For $\xi=0$, consists of $60$ separate segments, each an arc of a great circle on the unit sphere (each beginning and ending at the fixed point of our superflow).
\item[v)]For $0<\xi<\frac{\phi^3}{27}=\frac{2+\sqrt{5}}{27}$, consists of $20$ disjoint curves, each homeomorphic to a circle. (See Figure \ref{figure-pos}).
\item[vi)]For $\xi=\frac{2+\sqrt{5}}{27}$, is equal to $20$ isolated points.
\item[vii)]For $\xi>\frac{2+\sqrt{5}}{27}$, is empty.
\end{itemize}
\end{prop}

The surface $\mathbf{T}_{\xi}$ is not compact. To get a better visualization, we may consider the compact surface 
\begin{eqnarray*}
\mathbf{L}_{\xi}=\{\m{x}\in\mathbb{R}^{3}:\mathscr{V}(\m{x})+(x^2+y^2+z^2)^3=1+\xi\}.
\end{eqnarray*}
For $\xi>-1$ it is compact and connected, and obviously, $\mathbf{L}_{\xi}\cap\mathbf{S}=\mathbf{T}_{\xi}\cap\mathbf{S}$. Two cases are shown in Figures \ref{figure-neg} and \ref{figure-pos}.\\
 
Now  we are ready to give the full list of  points where the vector field vanishes. At these points $V(x,y,z)=x$. As it appears, these are precisely $62$ points: the endpoints of archs defined in the item iv), $30$ points in total; $20$ centres of all triangles; $12$ centres of all pentagons. And indeed, we can check directly that $\varpi\bl\varrho\bl\sigma=0\bl0\bl 0 $ for $(x,y,z)$ equal to:
\begin{itemize}
\item[P)] $12$ points  $(\pm\phi,\pm 1,0)$, $(0,\pm\phi,\pm 1)$, $(\pm 1,0,\pm\phi)$ (centres of pentagons).
\item[T)] $20$ points $(\pm 1,\pm 1,\pm 1)$, $(\pm\frac{1}{\phi},\pm\phi,0)$, $(\pm\phi,0,\pm\frac{1}{\phi})$, 
$(0,\pm\frac{1}{\phi},\pm\phi)$ (centres of triangles).
\item[E)] $24$ points $(\pm\phi,\pm\frac{1}{\phi},\pm 1)$ (this and all cyclic permutations) and $6$ points $(\pm 1,0,0)$, $(0,\pm 1,0)$, $(0,0,\pm 1)$ (all vertices). 
\end{itemize} 
\begin{Note} For the octahedral superflow \cite{alkauskas-super1}, we saw that the number of fixed points of the superflow on the unit sphere is equal to $26$, exactly the total sum of faces ($8$), vertices ($6$), and edges ($12$) of an octahedron. Here we encounter precisely the same situation - the total number of fixed points for the icosahedral superflow is exactly the sum of faces ($20$), vertices ($12$), and edges ($30$) of an icosahedron. We can say more - the index of the vector field at a fixed point is $+1$, if it corresponds to a face or a vertex, and equal to $-1$, if it corresponds to an edge, thus making the sum of all indices equal to the Euler characteristics of the sphere. 
\end{Note}
\subsection{Stereographic projection} To visualize the vector field $\varpi\bl\varrho\bl\sigma$ on the unit sphere, we use a stereograpic projection from the point $\mathbf{p}=(0,0,1)$ (where the vector field vanishes) to the plane $z=-1$. There is a canonical way to do a projection of a vector field as follows. Let $F(\m{x},t)=t^{-1}\phi(\m{x}t)$ be our flow on the unit sphere, and $\tau$ - the stereographic projection from the point $\mathbf{p}$ to the plane $z=-1$. We call a $2$-dimenional vector field $\mathcal{Y}$ in the plane $z=-1$ a \emph{stereographic projection} of the vector field $\mathcal{X}=\varpi\bl\varrho\bl\sigma$, if for a flow $G$,  associated with the vector field $\mathcal{Y}$, one has
\begin{eqnarray}
G(\tau(\m{x}),t)=\tau(F(\m{x},t))\text{ for all }\m{x}\text{ and }t\text{ small enough}.
\label{proj}
\end{eqnarray}
Let $\mathbf{a}=(x,y,z)$, $\mathbf{a}\neq\mathbf{p}$, $x^2+y^2+z^2=1$. The line $\overline{\mathbf{pa}}$ intersects the plane $z=-1$ at the point $\mathbf{b}$, which we call $\tau(\mathbf{a})$; this is, of course, a standard stereographic projection. Now, consider the line $(x+\varpi t,y+\varrho t, z+\sigma t)$, $t\in\mathbb{R}$. Via a stereographic projection (only one point of it belongs the unit sphere, but nevertheless stereographic projection works in the same way), this line is mapped into a smooth curve (in fact, a quadratic) $\mathcal{Q}(\mathbf{a})$ contained in the plane $z=-1$. Its tangent vector at the point $\mathbf{a}$ is a direction of a new vector in the plane $z=-1$.\\

The point $\mathbf{a}=(x,y,z)$ is mapped into a point $\tau(\mathbf{a})=(\alpha,\beta)=(\frac{2x}{1-z},\frac{2y}{1-z})$ (where $z=-1$), and the point $(x+\varpi t,y+\varrho t, z+\sigma t)$ - into a point
\begin{eqnarray*}
\alpha(t)\bl\beta(t)=\tau(t,\mathbf{a})=\frac{2x+2t\varpi}{1-z-t\sigma}\bl\frac{2y+2t\varrho}{1-z-t\sigma}.
\end{eqnarray*}
If $\gamma=\gamma_{x,y,z}$ and $\delta=\delta_{x,y,z}$ are the corresponding M\"{o}bius transformations on the right (as a functions in $t$, where $x,y,z$ being fixed, and consequently $\varpi,\varrho,\sigma$ also), so that $\alpha(t)=\gamma(t)$, $\beta(t)=\delta(t)$, then the equation for the curve $\mathcal{Q}(\mathbf{a})$ can be written as $\gamma^{-1}(\alpha)=\delta^{-1}(\beta)$, which shows immediately that it is a quadratic. Of course, we are interested only at its germ at $\tau(\mathbf{a})$. 
The inverse of $\tau$ is given by 
\begin{eqnarray}
\tau^{-1}(\alpha,\beta)=\frac{4\alpha}{\alpha^2+\beta^2+4}
\bl\frac{4\beta}{\alpha^2+\beta^2+4}\bl\frac{\alpha^2+\beta^2-4}{\alpha^2+\beta^2+4}.\label{ab-inv}
\end{eqnarray}
Now,
\begin{eqnarray}
\frac{\d}{\d t}\tau(t,\mathbf{a})\Big{|}_{t=0}=\frac{2\varpi-2\varpi z+2x\sigma}{(1-z)^2}\bl
\frac{2\varrho-2\varrho z+2y\sigma}{(1-z)^2}.
\label{true-vec}
\end{eqnarray}
If we use now a substitution (\ref{ab-inv}), that is, plug $(x,y,z)\mapsto \tau^{-1}(\alpha,\beta)$ into the above, we obtain exactly the stereographic projection of the vector field $\varpi\bl\varrho\bl\sigma$. Indeed, $F(\m{x},t)=x+t\varpi+O(t^2)\bl y+t\varrho+O(t^2)\bl z+t\sigma+O(t^2)$. So, $\tau(F(\m{x},t))=\tau(t,\mathbf{a})+O(t^2)$, and the answer follows.\\

To get a better visualization, we multiply (\ref{true-vec}) by $1-z$; this does not change nor orbits neither directions, but changes the flow, of course. We thus arrive at the vector field
\begin{eqnarray}
2\varpi+\frac{2x}{1-z}\sigma\bl 2\varrho+\frac{2y}{1-z}\sigma. 
\label{exprr}
\end{eqnarray}
Or, minding (\ref{ab-inv}), in terms of $(\alpha,\beta)$ this reads as
\begin{eqnarray*}
2\varpi\Big{(}\tau^{-1}(\alpha,\beta)\Big{)}+\alpha\sigma\Big{(}\tau^{-1}(\alpha,\beta)\Big{)}
\bl 2\varrho\Big{(}\tau^{-1}(\alpha,\beta)\Big{)}+\beta\sigma\Big{(}\tau^{-1}(\alpha,\beta)\Big{)}=\Pi(\alpha,\beta)\bl\Theta(\alpha,\beta).
\end{eqnarray*}
This is a vector field on the $(\alpha,\beta)$ plane. This is inessential for visualisation, but for purposes described in Subsection \ref{orth-sub}, note that the true projection of the vector field $\mathcal{X}$ is 
\begin{eqnarray}
\frac{1}{8}(\alpha^2+\beta^2+4)\Pi(\alpha,\beta)\bl \frac{1}{8}(\alpha^2+\beta^2+4)\Theta(\alpha,\beta), 
\label{ico-pro}
\end{eqnarray}
and is given by a pair of rational functions, but not $2-$ homogenic anymore.\\

 If $\varpi\bl\varrho\bl\sigma$ vanishes on the unit sphere at a point $\mathbf{a}$, so does $\Pi\bl\Theta$ at a point $\tau(\mathbf{a})$. In the other direction - suppose $\Pi=0$, $\Theta=0$ for a certain $(\alpha,\beta)$. Then there is a point $(x,y,z)\neq (0,0,1)$ on the unit sphere such that both expressions in (\ref{exprr}) vanish. Recall that we also have $x\varpi+y\varrho+z\sigma=0$. Now, the determinant of the  matrix
\begin{eqnarray*}
\begin{pmatrix}
1-z & 0 & x\\
0 & 1-z & y\\
x & y & z 
\end{pmatrix}
\end{eqnarray*}
is equal to $-(z-1)^2\neq 0$. This implies $\varpi\bl\varrho\bl\sigma=0\bl 0\bl 0$.\\

\begin{figure}
\epsfig{file=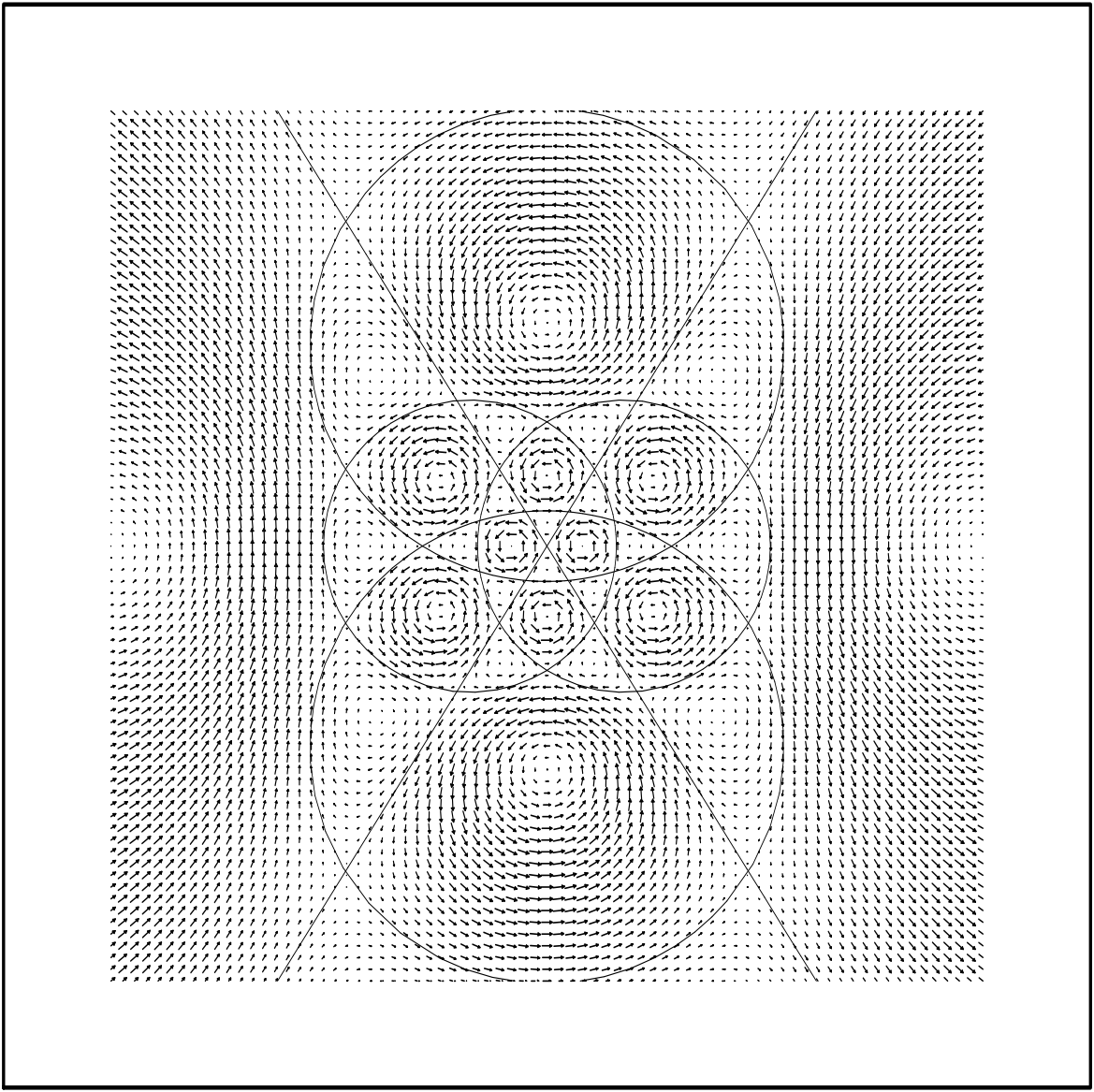,width=440pt,height=440pt,angle=-90}
\caption{Stereographic projection of the icosahedral superflow, $-7\leq\alpha,\beta\leq 7$. Orbits inside each pentagon or triangle are non-singular (except for a center which is a fixed point and an orbit on its own), and this case is treated in Section \ref{sec-non-sing}. Further, each segment is a separate orbit; this case is simpler and is treated in Section \ref{sec-sing}. }
\label{figure2}
\end{figure}

\begin{figure}
\epsfig{file=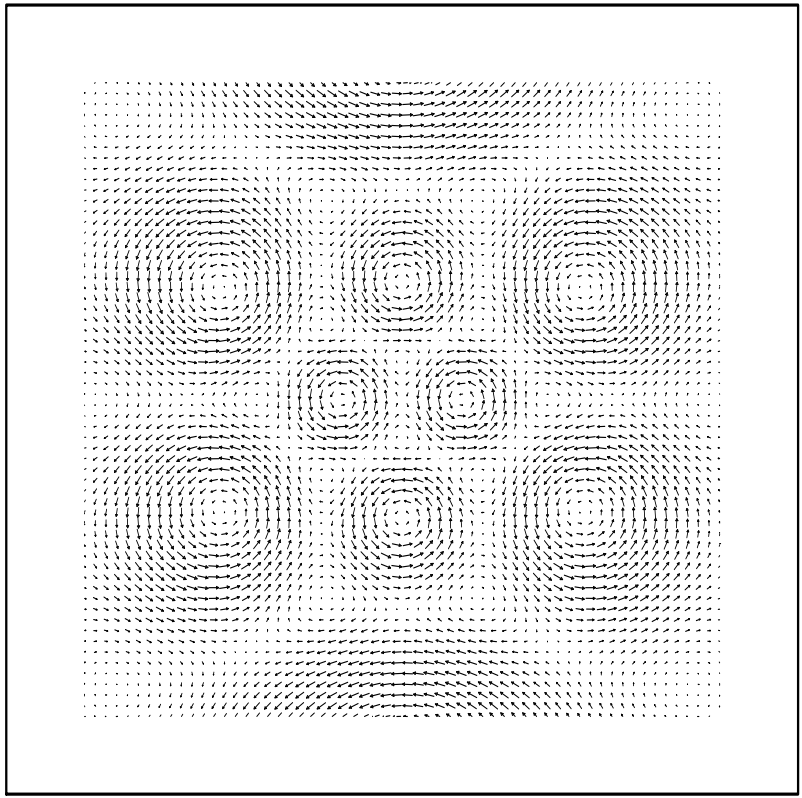,width=300pt,height=300pt,angle=-90}
\caption{A close-up of the stereographic projection of the icosahedral superflow, $-3\leq\alpha,\beta\leq 3$.}.
\label{figure3}
\end{figure}

The intersections of the plane $y=\phi x$ and the plane $y=-\phi x$ with the unit sphere are mapped by $\tau$ to the lines $\alpha=\phi\beta$ and $\alpha=-\phi\beta$, respectively. The intersection of the plane $z=\pm\phi y$ and the unit sphere is parametrized by $\Big{(}\sin\vartheta,\frac{\cos\vartheta}{\sqrt{\phi^2+1}},\pm\frac{\phi\cos\theta}{\sqrt{\phi^2+1}}\Big{)}$, $\vartheta\in[0,2\pi]$. So, in the $(\alpha,\beta)$ plane we get a parametrization
\begin{eqnarray*}
\frac{2\sqrt{\phi^2+1}\sin\vartheta}{\sqrt{\phi^2+1}\mp\phi\cos\theta}\bl
\frac{2\cos\vartheta}{\sqrt{\phi^2+1}\mp\phi\cos\theta},\quad\vartheta\in[0,2\pi].
\end{eqnarray*}
This is a circle with center $(0,2\phi)$ and a radius $2\sqrt{\phi^2+1}$ for an upper sign, and a circle with the center $(0,-2\phi)$ with the same radius for a lower sign, respectively. Analogously we find that the intersection of two planes $x=\pm\phi z$ with the unit sphere are mapped by $\tau$ to two circles with radius $\phi^{-1}\sqrt{\phi^2+1}$ and centres $(\pm2\phi^{-1},0)$, respectively. The Figures \ref{figure2} and \ref{figure3} demonstrates the vector field and six boundary curves in the $(\alpha,\beta)$ plane.\\

\begin{figure}
\epsfig{file=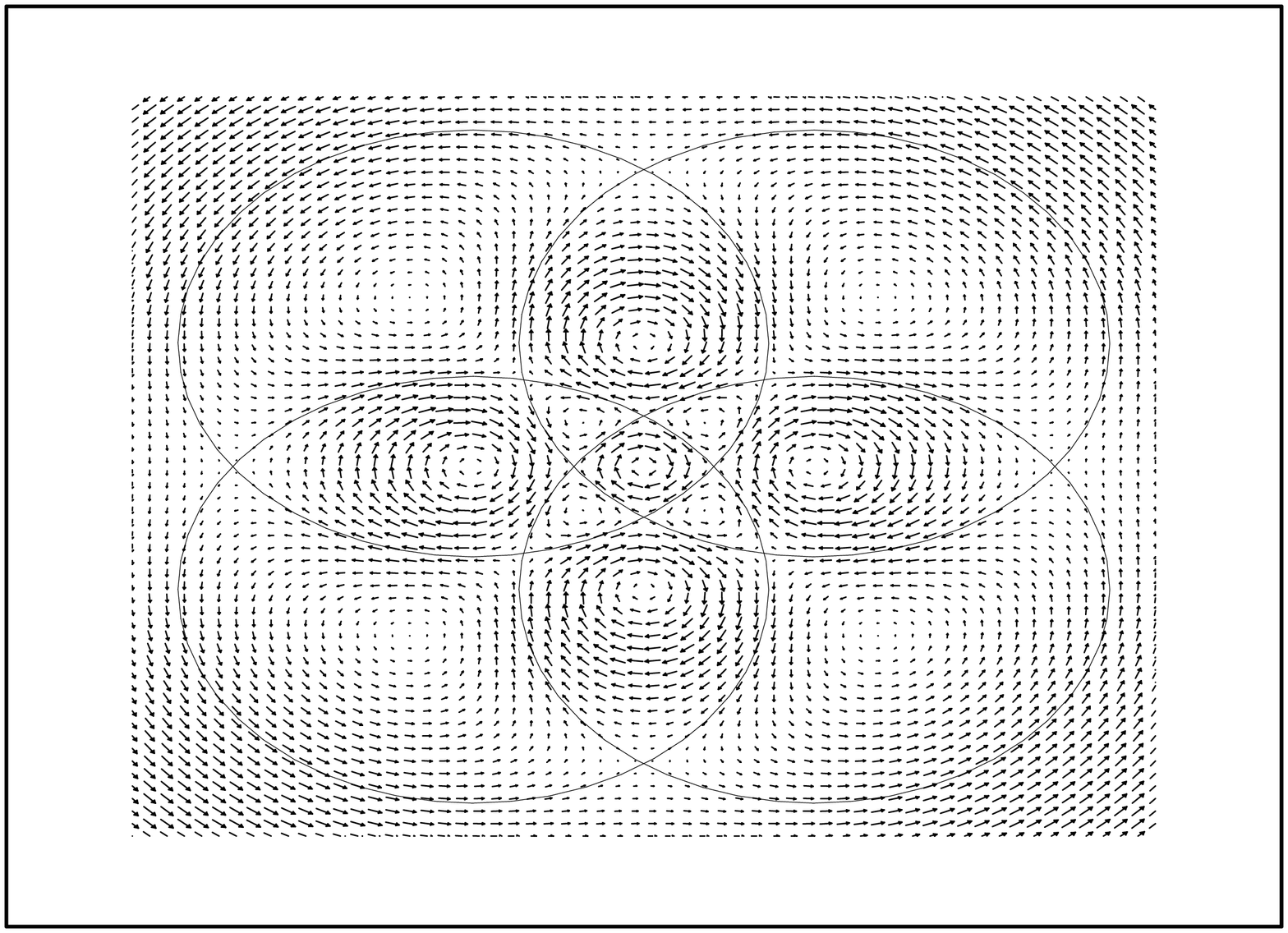,width=300pt,height=300pt,angle=-90}
\caption{Stereographic projection of the octahedral superflow, $-3\leq\alpha,\beta\leq 3$.}
\label{figure4}
\end{figure}

We further present a stereographic projection of the octahedral superflow, treated in detail in \cite{alkauskas-super1} (see Theorem \ref{thm1}). This time the singular orbits are curves $\{x^2+y^2+z^2=1,x^4+y^4+z^4=\frac{1}{2}\}$.  This is a reducible case, since
\begin{eqnarray*}
(x+y+z)(x+y-z)(x-y+z)(x-y-z)=2(x^4+y^4+z^4)-(x^2+z^2+z^2)^2=0.
\end{eqnarray*}
So, these singular orbits are intersections of the unit sphere with $4$ planes. A direct calculation shows that stereographically these four circles in the $(\alpha,\beta)$ plane map into $4$ circles with radii $\sqrt{12}$, and centres $(\pm 2,\pm 2)$ (signs are independent). Figure \ref{figure4} shows the setup.

\subsection{Orthogonal and other birational projections}
\label{orth-sub}
In the same vein, as another, much simpler example how to visualize a $3$-dimensional vector field, we will give an orthogonal projection of the vector field
\begin{eqnarray*}
\mathcal{X}=\varpi\bl\varrho\bl\sigma=yz\bl xz\bl xy,
\end{eqnarray*}
which gives rise to the tetrahedral superflow $\phi_{\widehat{\mathbb{T}}}$ , as described by Theorem \ref{thm1}, part i.\\

Consider the surface $\mathscr{Q}=\{(x,y,z)\in\mathbb{R}^{3}:x^2-y^2=1,x>0\}$, which is one of a two connected components of a quadratic. Since $p^2-q^2$ is the first integral of the differential system (\ref{sys-in}), the vector field $\mathcal{X}$ is tangent to $\mathscr{Q}$. We will project now orthogonally $\mathcal{X}$ to the plane $x=0$, exactly as in the previous subsection. Let $\tau(\m{x})=\tau(x,y,z)=(y,z)$ be this projection. Let $F(\m{x},t)=t^{-1}\phi(\m{x}t)$ be a flow on $\mathscr{Q}$. Equally, we call a $2$-dimensional vector field $\mathcal{Y}$ in the plane $x=0$ an \emph{orthogonal projection} of the vector field $\mathcal{X}$, if for a flow $G$, associated with the vector field $\mathcal{Y}$, one has (\ref{proj}).   \\

So, let $\mathbf{a}=(x,y,z)$ belongs to $\mathscr{Q}$. It is mapped to the point $\tau(\m{a})=(\alpha,\beta)=(y,z)$ (where $x=0$), and the point $(x+\varpi t,y+\varrho t,z+\sigma t)$ on the tangent line - to a point $\tau(t,\m{a})=(y+\varrho t,z+\sigma t)$. The inverse of $\tau$ is given by
\begin{eqnarray*}
\tau^{-1}(\alpha,\beta)=\sqrt{1+\alpha^2}\bl\alpha\bl\beta.
\end{eqnarray*}
Now, obviously, 
\begin{eqnarray*}
\frac{\d}{\d t}\tau(t,\m{a})\Big{|}_{t=0}=\varrho\bl\sigma.
\end{eqnarray*}   
As before, we now plug $\tau^{-1}(\alpha,\beta)\mapsto(x,y,z)$ into the above to obtain an orthogonal projection of our vector field $\mathcal{X}$. Thus, we have the vector field
\begin{eqnarray}
\Pi(\alpha,\beta)\bl\Theta(\alpha,\beta)=\beta\sqrt{1+\alpha^2}\bl \alpha\sqrt{1+\alpha^2}.
\label{pi-theta}
\end{eqnarray} 
On a surface $\mathscr{Q}$, the orbits of the flow with the vector field $\mathcal{X}$ are intersections of $\mathscr{Q}$ with $x^2-z^2=\xi$, $\xi\in\mathbb{R}$. Generically, these  are elliptic curves, except for cases $\xi=0$ and $\xi=1$. These exceptions correspond to intersection of $\mathscr{Q}$ with the planes $x=\pm z$ and $z=\pm y$.  In the $(\alpha,\beta)$ plane, this corresponds to a quadratic $\beta^2-\alpha^2=1$, and two lines $\alpha=\pm \beta$, respectively. In general, the orbits of the flow with the vector field $\Pi\bl \Theta$ are quadratics $\alpha^2-\beta^2=\phi$, $\phi\in\mathbb{R}$. Thus, we get a standard partition of the plane into hyperbolas having two fixed lines as their asymptotes. \\

Note that, differently from a previous subsection and the vector field (\ref{ico-pro}), the vector field $\Pi\bl\Theta$ in this subsection is not rational, since the map $\tau$ is not a birational map. However, the projected flow can be described in rational terms, too, if we use a projection to the plane $x-y=0$ rather that to the plane $x=0$. This time not only $\mathcal{Q}$ is projected injectively, but the whole hyperbolic cylinder $\widehat{\mathcal{Q}}=\{(x,y,z)\in\mathbb{R}^{3}:x^2-y^2=1\}$.\\

 Indeed, suppose, a point $(\alpha,\alpha,\beta)$ in the $x-y=0$ plane  is described by global coordinates $\alpha\bl\beta$. Thus, if $(x,y,z)\in\widehat{\mathcal{Q}}$, then the orthogonal projection to the plane $x-y=0$, and its inverse, are given by
\begin{eqnarray*}
\hat{\tau}(x,y,z)=\frac{x+y}{2}\bl z,\quad
\hat{\tau}^{-1}(\alpha,\beta)=\Big{(}\alpha+\frac{1}{4\alpha}\Big{)}\bl \Big{(}\alpha-\frac{1}{4\alpha}\Big{)}\bl\beta.
\end{eqnarray*}
In the same fashion as before, in the $(\alpha, \beta)$ plane we get a vector field
\begin{eqnarray}
\hat{\Pi}(\alpha,\beta)\bl\hat{\Theta}(\alpha,\beta)=\alpha\beta\bl\alpha^{2}-\frac{1}{16\alpha^2}.
\label{al-be}
\end{eqnarray}
Thus, using the result in \cite{alkauskas-super1}, we get
\begin{prop}
\label{jac-alt}
The $2$-dimensional vector fields $\Pi\bl\Theta$ and $\hat{\Pi}\bl\hat{\Theta}$, as given by (\ref{pi-theta}) and (\ref{al-be}), which produce flows, but not projective flows any longer, can be explicitly integrated in terms of Jacobi elliptic functions.  
\end{prop}
However, differently from $\Pi\bl\Theta$, $\hat{\Pi}\bl\hat{\Theta}$ is given by rational functions. So, though the vector field $\hat{\Pi}\bl\hat{\Theta}$ is not given by a pair of $2-$homogenic functions, it still can be described via a projective flow $\phi_{\hat{\mathbb{T}}}$ defined in Theorem \ref{thm1}.\\

Further, let $x^2-y^2=1$, $x^2-z^2=\xi$. Thus, we take any point on one particular orbit.  Let, as before, $\alpha\bl\beta=\frac{x+y}{2}\bl z$. We know that $\alpha+\frac{1}{4\alpha}=x$. This can be checked directly, using the relation $\frac{1}{x+y}=x-y$. Thus, we obtain the equation
\begin{eqnarray*}
\mathscr{W}(\alpha,\beta)=\Big{(}\alpha+\frac{1}{4\alpha}\Big{)}^2-\beta^2=\xi.
\end{eqnarray*} 
This is the equation for the generic orbit for the flow with the vector field $\hat{\Pi}\bl\hat{\Theta}$, and for $\xi\neq 0,1$ it is an elliptic curve. We can double verify that the following holds:
\begin{eqnarray*}
\mathscr{W}_{\alpha}\hat{\Pi}+\mathscr{W}_{\beta}\hat{\Theta}\equiv 0.
\end{eqnarray*}

To this account that we make one elementary observation. Let
\begin{eqnarray*}
\sum\limits_{k=1}^{n}f_{k}\frac{\p}{\p x_{k}}
\end{eqnarray*}
be any $n$-dimensional vector field. Let us define
\begin{eqnarray*}
\tilde{f}_{k}(x_{1},x_{2},\ldots,x_{n},z)=
z^{2}f_{k}\Big{(}\frac{x_{1}}{z},\frac{x_{2}}{z},\ldots,\frac{x_{n}}{z}\Big{)},\quad k=1,\ldots,n.
\end{eqnarray*}
Then these functions are $2$-homogenic, and so the vector field
\begin{eqnarray}
\sum\limits_{k=1}^{n}\tilde{f}_{k}\frac{\p}{\p x_{k}}+0 \frac{\p}{\p z}
\label{confine}
\end{eqnarray}
gives rise to an $(n+1)$-dimensional projective flow. Note that $z=1$ is its integral surface. If we confine the projective flow (\ref{confine}) to the plane $z=1$, we recover the initial flow. Hence, we have
\begin{prop}
\label{prop-ner}
Any flow in $\mathbb{R}^{n}$ is a section on a hyperplane of a $(n+1)$-dimensional projective flow. In other words, for any flow $F(\m{x},t):\mathbb{R}^{n}\times\mathbb{R}\mapsto\mathbb{R}^{n}$, there exists a projective flow $\phi:\mathbb{R}^{n+1}\mapsto\mathbb{R}^{n+1}$, such that
\begin{eqnarray*}
F(\m{x},t)=\frac{1}{t}\phi\Big{(}(\m{x},1)t\Big{)}=\frac{1}{t}\phi\big{(}\m{x}t,t\big{)}.
\end{eqnarray*}
\end{prop} 
This holds even if the flow $F$ is defined for $\m{x}$ belonging to an open region of $\mathbb{R}^{n}$, and $t$ small enough (depending on $\m{x}$). On the other hand, the projective flow (\ref{confine}) is of a special kind, and hence the variety of projective flows in dimension $(n+1)$ is bigger than that of general flows in dimension $n$. So, philosophically, we may consider (\ref{funk}) for $\m{x}\in\mathbb{R}^{n}$ as a general flow equation in dimension $n-\frac{1}{2}$. \\

However, the scenario in Proposition \ref{prop-ner} is not the only way to obtain a flow as a projection of a projective flow, as we have seen.\\

Suppose $\varpi\bl\varrho\bl\sigma$ is a triple of $2$-homogenic rational functions, and let $\mathscr{W}$ be a homogenic function and the first integral of the flow, as in (\ref{first-int}). Suppose now, an integral surface $\mathcal{S}$, given by $\mathscr{W}(x,y,z)=1$, is of genus $0$. For example, in Proposition \ref{prop-ner} this is just $z=\mathrm{const.}$ \\

 In the most general case, by a \emph{projection of a vector field} we mean the following. Consider a rational parametrization $\pi:\mathbb{R}^{2}\mapsto\mathcal{S}$, $\pi:(\alpha,\beta)\mapsto\big{(}f(\alpha,\beta),g(\alpha,\beta),h(\alpha,\beta)\big{)}$. Suppose that its inverse $\pi^{-1}:\mathcal{S}\mapsto\mathbb{R}^2$ is also rational, and is given by rational functions $\pi^{-1}:(x,y,z)\mapsto(A(x,y,z),B(x,y,z))$. Thus, the projection of our vector field is a new vector field in the $(\alpha,\beta)$ plane given by
\begin{eqnarray*}
\Big{(}A_{x}\varpi+A_{y}\varrho+A_{z}\sigma, 
B_{x}\varpi+B_{y}\varrho+B_{z}\sigma\Big{)}\Big{|}_{\big{(}(f(\alpha,\beta),g(\alpha,\beta),h(\alpha,\beta)\big{)}\mapsto(x,y,z)}.
\end{eqnarray*}
In both cases dealt with in this paper (a sphere in icosahedral and octahedral superflow cases, and a hyperbolic cylinder in a tetrahedral superflow case, orthogonal projection to the plane $x-y=0$), surfaces are of genus $0$, and are birationally parametrized. 
\section{The differential system}
\subsection{Algebraic curve parametrization}
To integrate the vector field (\ref{vec-field}) in closed terms, we need to solve the following differential system:
\begin{eqnarray}\scriptsize
\left\{\begin{array}{l}
p'=-(5-\sqrt{5})qr^5-(5+\sqrt{5})q^5r+20q^3r^3-(10+10\sqrt{5})p^2qr^3-(10-10\sqrt{5})p^2q^3r+10p^4qr,\\
q'=-(5-\sqrt{5})rp^5-(5+\sqrt{5})r^5p+20r^3p^3-(10+10\sqrt{5})q^2rp^3-(10-10\sqrt{5})q^2r^3p+10q^4rp,\\
r'=-(5-\sqrt{5})pq^5-(5+\sqrt{5})p^5q+20p^3q^3-(10+10\sqrt{5})r^2pq^3-(10-10\sqrt{5})r^2p^3q+10r^4pq,\\
p^2+q^2+r^2=1,\quad (\phi^2 p^2-q^2)(\phi^2 q^2-r^2)(\phi^2 r^2-p^2)=\xi.
\end{array}
\right.
\label{sys1}
\end{eqnarray}\normalsize
Note that in first three equations no denominators occur, since $p^2+q^2+r^2=1$.\\

Now we will make few steps in investigating the system (\ref{sys1}). Of course, since $\mathscr{W}$ and $\mathscr{V}$ are the first integrals, the second and the third equality follows, as always in (\ref{sys-in}), from the first and the last two. Let $\phi^2p^2-q^2=P$, $\phi^2q^2-r^2=Q$, $\phi^2r^2-p^2=R$. The inverse is given by
\begin{eqnarray}
\left\{\begin{array}{l}
4\phi p^2=\phi^2 P+Q+\phi^{-2}R,\\
4\phi q^2=\phi^{-2}P+\phi^2Q+R,\\
4\phi r^2=P+\phi^{-2}Q+\phi^2R. 
\end{array}
\label{aha}
\right.
 \end{eqnarray}
The last two identities of the system (\ref{sys1}) in terms of $P,Q,R$ can be rewritten as
\begin{eqnarray}
P+Q+R=\phi,\quad PQR=\xi.
\label{rys}
\end{eqnarray}

From the first and the second identities of (\ref{sys1}) we imply, by a direct calculation,
\begin{eqnarray*}
2\phi^2pp'-2qq'=(40-8\sqrt{5})pqr(\phi^2 p^2-q^2)\Big{(}\phi^{2}q^2+p^2
-(\phi^{2}+1)r^2\Big{)}.
\end{eqnarray*}
Since $P'=2\phi^2pp'-2qq'$, this gives
\begin{eqnarray*}
P'^2=1280\phi^{-2}\cdot p^2q^2r^2\cdot P^2\cdot(Q-R)^2.
\end{eqnarray*}
From (\ref{rys}) we get $(Q-R)^2=(Q+R)^2-4QR=(\phi-P)^2-\frac{4\xi}{P}$, and so
\begin{eqnarray}
\frac{P'^{2}}{P^4-2\phi P^3+\phi^2 P^2-4\xi P}=1280\phi^{-2}\cdot p^2q^2r^2.
\label{recu}
\end{eqnarray}
Let us now put $\alpha=P^2Q+Q^2R+R^2P$, $\beta=PQ^2+QR^2+RP^2$.
Now (\ref{aha}), after multiplying all three identities throughout, gives
\begin{eqnarray}
64\phi^{3}p^2q^2r^2&=&22PQR+P^3+Q^3+R^3+(6+\phi)\alpha+(7-\phi)\beta\nonumber\\
&=&22PQR+P^3+Q^3+R^3+\frac{13}{2}(\alpha+\beta)+\frac{\sqrt{5}}{2}(\alpha-\beta).
\label{bas}
\end{eqnarray}
We will express the right hand side of (\ref{bas}) in terms of $P$ using (\ref{rys}). First,
\begin{eqnarray}
P^3+Q^3+R^3&=&P^3+(Q+R)^3-3QR(Q+R)=P^3+(\phi-P)^3-\frac{3\xi(\phi-P)}{P}\nonumber\\
&=&\phi^3-3\phi^2P+3\phi P^2-\frac{3\xi\phi}{P}+3\xi,\label{pirms}
\end{eqnarray}
and
\begin{eqnarray}
(P+Q+R)^3&=&P^3+Q^3+R^3+6PQR+3(\alpha+\beta)\Rightarrow\nonumber\\
\phi^3&=&\phi^3-3\phi^2P+3\phi P^2-\frac{3\xi\phi}{P}+3\xi+6\xi+3(\alpha+\beta)\nonumber\\
\alpha+\beta&=&\phi^2P-\phi P^2+\frac{\xi\phi}{P}-3\xi.\label{antrs}
\end{eqnarray}
Now, let us plug (\ref{pirms}) and (\ref{antrs}) into (\ref{bas}). We obtain:
\begin{eqnarray*}
128\phi^{3}p^2q^2r^2&=&44\xi+2(\phi^3-3\phi^2P+3\phi P^2-\frac{3\xi\phi}{P}+3\xi)\\
&+&13(\phi^2P-\phi P^2+\frac{\xi\phi}{P}-3\xi)+\sqrt{5}(\alpha-\beta)\\
&=&11\xi+2\phi^3+7\phi^2P-7\phi P^2+\frac{7\xi\phi}{P}+\sqrt{5}(\alpha-\beta).
\end{eqnarray*}
(We can double-verify the validity of the last formula with MAPLE by plugging $\xi\mapsto PQR$, $\phi\mapsto P+Q+R$, thus making the expression on the right $3-$homogenic in $P,Q,R$, and then using (\ref{aha}) to see that this indeed is equal to $128\phi^3p^2q^2r^2$).\\

Plugging now this into (\ref{recu}), we obtain
\begin{eqnarray*}
\frac{P'^2\phi^{5}}{10(P^4-2\phi P^3+\phi^2P^2-4\xi P)}=
11\xi+2\phi^3+7\phi^2P-7\phi P^2+\frac{7\xi\phi}{P}+\sqrt{5}(\alpha-\beta).
\end{eqnarray*}
Or, more conveniently,
\begin{eqnarray*}
\frac{P'^2\phi^{5}}{10(P^4-2\phi P^3+\phi^2P^2-4\xi P)}-11\xi-2\phi^3-7\phi^2P+7\phi P^2-\frac{7\xi\phi}{P}=\sqrt{5}(\alpha-\beta).
\end{eqnarray*}
The square of the right hand side is a rational function in $P$, and this gives us a polynomial equation for a pair $(X,Y)=(P',P)$. Indeed, $(\alpha-\beta)^2=(P-Q)^2(Q-R)^2(R-P)^2$ is the discriminant of the polynomial
\begin{eqnarray*}
(\mathbf{x}-P)(\mathbf{x}-Q)(\mathbf{x}-R)=(\mathbf{x}-P)\Big{(}\mathbf{x}^2-(\phi-P)\mathbf{x}+\frac{\xi}{P}\Big{)},
\end{eqnarray*}
and so 
\begin{eqnarray*}
(\alpha-\beta)^2=\Big{(}2P^2-\phi P+\frac{\xi}{P}\Big{)}^2\Big{(}P^2-2\phi P+\phi^2-\frac{4\xi}{P}\Big{)}.
\end{eqnarray*}
Thus, finally, if $(X,Y)=(P',P)$, then
\begin{eqnarray*}
& &\Bigg{(}\frac{Y^2\phi^{5}}{10(X^4-2\phi X^3+\phi^2X^2-4\xi X)}-11\xi-2\phi^3-7\phi^2X+7\phi X^2-\frac{7\xi\phi}{X}\Bigg{)}^2\\
&=&
5\Big{(}2X^2-\phi X+\frac{\xi}{X}\Big{)}^2\Big{(}X^2-2\phi X+\phi^2-\frac{4\xi}{X}\Big{)}.
\end{eqnarray*}
This can be rewritten as
\begin{eqnarray}
& &X\Bigg{(}Y^2\phi^{5}+10\Big{(}7\phi X^3-7\phi^2X^2-(11\xi+2\phi^3)X-7\xi\phi\Big{)}\cdot\Big{(}X^3-2\phi X^2+\phi^2X-4\xi\Big{)}\Bigg{)}^2\nonumber\\
&=&
500\Big{(}2X^3-\phi X^2+\xi\Big{)}^2\Big{(}X^3-2\phi X^2+\phi^2X-4\xi\Big{)}^{3}.
\label{comp}
\end{eqnarray}
Note that all the steps are valid for $Q,R$ in place of $P$. Thus, we have proved the following.
\begin{prop}
\label{prop2}
Let us define the polynomials
\begin{eqnarray*}
f_{\xi}(X)&=&7\phi X^3-7\phi^2X^2-(11\xi+2\phi^3)X-7\xi\phi,\\
g_{\xi}(X)&=&X^3-2\phi X^2+\phi^2X-4\xi,\\
h_{\xi}(X)&=&2X^3-\phi X^2+\xi,
\end{eqnarray*}
where $\xi\in\mathbb{R}$, $-\frac{2+\sqrt{5}}{5}<\xi<\frac{2+\sqrt{5}}{27}$, is arbitrary but fixed. Let $(X,Y)=(P,P')$, $(Q,Q')$, or $(R,R')$ (a function and its derivative). Then a pair of functions $(X,Y)$ parametrizes the algebraic curve
\begin{eqnarray}
X\Big{(}Y^2\phi^{5}+10f_{\xi}(X)g_{\xi}(X)\Big{)}^2=500\cdot h_{\xi}^{2}(X)\cdot g^{3}_{\xi}(X).
\label{comp2}
\end{eqnarray}
Consider the above as a fourth degree polynomial in $Y$. Then it is irreducible in $\mathbb{Q}(\xi,X)[Y]$.\\

For $\xi$ in the range we consider, $\xi\neq 0$, (\ref{comp2}) defienes a curve of arithmetic genus $7$.
\end{prop}
Compare this to (\cite{alkauskas-super1}, Proposition 7), valid in the context of the second superflow $\phi_{\mathbb{O}}$ described by Theorem \ref{thm1}.
\begin{proof}The proof is presented above, but we can symbolically double-verify the claim of the Proposition with the help of MAPLE. Let us make both sides of (\ref{comp}) homogeneous $30$th degree polynomials in $p,q,r$ as follows. Let us put $X\mapsto P$, $\phi\mapsto P+Q+R$ (everywhere except for a factor $\phi^5$ just after $Y^2$), $\xi\mapsto PQR$, thus making both expressions on the right and the left homogenic polynomials in $P,Q,R$ of degree $15$; except for a moment we pay no attention to $Y$, and its factor $\phi^5$ remains intact. Now, let us further use a substitution $P\mapsto\phi^2p^2-q^2$ (this time $\phi$ is a constant), $Q\mapsto\phi^2q^2-r^2$, $R\mapsto\phi^2q^2-p^2$. Moreover, let $Y\mapsto2\phi^2p p'-2qq'$, where $p'$ and $q'$ should be substituted by the $6$th degree polynomials, the first and the second entry of the differential system (\ref{sys1}), respectively. this makes $Y$ a $7$th degree homogeneous function in $p,q,r$, and both sides of the equation (\ref{comp}) - homogeneous functions in $p,q,r$ of degree $30$. MAPLE shows symbolically, that the difference of both sides is indeed $0$, thus proving again (\ref{comp}).

\end{proof}
\subsection{Reduction}Let $X=P$, $Q$ or $R$, and let us introduce a function $\Upsilon(t)$ by the identity 
\begin{eqnarray*}
\frac{X^3-\phi X^2-\xi}{X}=\Upsilon(t).
\end{eqnarray*} 
This trick is completely analogous to the one used in (\cite{alkauskas-super1}, Section 7) - if $P,Q,R$ are three distinct roots of the above, then $P+Q+R=\phi$, $PQR=\xi$, as desired. We will see that the function $\Upsilon$ satisfies a much simpler differential equation than the one given by (\ref{comp}).\\

By a direct calculation,
\begin{eqnarray*}
\Upsilon'=\frac{h_{\xi}(X)}{X^2}X'.
\end{eqnarray*}
Recall that $\xi\in\mathbb{R}$ is fixed. We claim that there exist polynomials $\mathfrak{p}_{\xi}(\mathbf{x})=\mathfrak{p}(\mathbf{x}),\mathfrak{q}_{\xi}(\mathbf{x})=\mathfrak{q}(\mathbf{x})\in\mathbb{R}[\mathbf{x}]$ of degrees $k$ and $\ell$, respectively, such that
\begin{eqnarray}
\big{(}\Upsilon'^2\phi^{5}+\mathfrak{p}(\Upsilon)\big{)}^2=\mathfrak{q}(\Upsilon).
\label{simpli}
\end{eqnarray}
In terms of $X$, this reads as 
\begin{eqnarray*}
\Bigg{(}\frac{h^2_{\xi}(X)}{X^4}X'^2\phi^{5}+\mathfrak{p}\Big{(}\frac{X^3-\phi X^2-\xi}{X}\Big{)}\Bigg{)}^2=\mathfrak{q}\Big{(}\frac{X^3-\phi X^2-\xi}{X}\Big{)}.
\end{eqnarray*}
Let, as before, $Y=X'$. This can be rewritten as
\begin{eqnarray}
\Bigg{(}X^{k-4}h^2_{\xi}(X)Y^2\phi^{5}+X^{k}\mathfrak{p}\Big{(}\frac{X^3-\phi X^2-\xi}{X}\Big{)}\Bigg{)}^2=X^{2k}\mathfrak{q}\Big{(}\frac{X^3-\phi X^2-\xi}{X}\Big{)}.
\label{uup}
\end{eqnarray}
According to Proposition \ref{prop2}, the $4$th degree polynomial (\ref{comp2}) in $Y$ is irreducible. So, in order (\ref{uup}) to hold, it should be a multiple of (\ref{comp2}). The degree of $X^{k}\mathfrak{p}(\ldots)$ in big brackets is equal to $3k$, while the degree of $X^{k-4}h_{\xi}^{2}(X)$ is $k+2$. The difference of these two numbers should be equal to $6$, the degree of $f_{\xi}(X)g_{\xi}(X)$. And so the only possible choice is $k=4$. Thus, (\ref{uup}) can be written as
\begin{eqnarray*}
X\Bigg{(}h^2_{\xi}(X)Y^2\phi^{5}+X^{4}\mathfrak{p}\Big{(}\frac{X^3-\phi X^2-\xi}{X}\Big{)}\Bigg{)}^2=X^{9}\mathfrak{q}\Big{(}\frac{X^3-\phi X^2-\xi}{X}\Big{)}.
\label{ups}
\end{eqnarray*}
Again, comparing the degrees to those of (\ref{comp2}), we see that $\ell=9$. So, we need to find polynomials $\mathfrak{p}$ of degree $4$ and $\mathfrak{q}$ of degree $9$ such that
\begin{eqnarray}
X^{4}\mathfrak{p}\Big{(}\frac{X^3-\phi X^2-\xi}{X}\Big{)}&=&10h_{\xi}^{2}(X)f_{\xi}(X)g_{\xi}(X),\label{math-p}\\
X^{9}\mathfrak{q}\Big{(}\frac{X^3-\phi X^2-\xi}{X}\Big{)}&=&
500h_{\xi}^{6}(X)g_{\xi}^{3}(X).\label{math-q}
\end{eqnarray} 
Then (\ref{simpli}) is proved. We will solve the first question with MAPLE, and will do the second question by hand due to a nice compatibility of both right sides; namely, that $(h^2g)^3=h^6g^{3}$. The first task is done recurrently: let us write $\mathfrak{p}(X)=aX^4+bX^3+cX^2+dX+e$, and  first compare the highest degrees of $X$ in (\ref{math-p}), then degrees one lower, and so on. We thus find the unique solution
\begin{eqnarray*}
\mathfrak{p}(X)=\Big{(}14\phi X-22\xi-4\phi^3\Big{)}
\Big{(}20X^3+5\phi^2 X^2-90\phi\xi X-135\xi^2-20\phi^3\xi\Big{)};
\end{eqnarray*}
(for convenience, we also factor into polynomials). Let $\mathfrak{l}(X)$ be the first linear polynomial, and $\mathfrak{t}(X)$ - the second cubic factor. We note that both factors can be viewed as homogenic, if $X,\xi,\phi$ are given weights $2,3,1$, respectively. This is what should be expected. We can verify now that
\begin{eqnarray*}
\mathfrak{q}(X)=4\mathfrak{t}^{3}(X).
\end{eqnarray*}  
Indeed, suppose that (\ref{math-p}) holds, and we will prove (\ref{math-q}). First, by a direct calculation,
\begin{eqnarray}
\mathfrak{l}\Big{(}\frac{X^3-\phi X^2-\xi}{X}\Big{)}=\frac{2f_{\xi}(X)}{X}.
\label{ties}
\end{eqnarray}
Thus,
\begin{eqnarray*}
X^{9}\mathfrak{q}\Big{(}\frac{X^3-\phi X^2-\xi}{X}\Big{)}
&&=4X^{9}\cdot\mathfrak{t}^{3}\Big{(}\frac{X^3-\phi X^2-\xi}{X}\Big{)}\\
&&=4X^{9}\frac{\mathfrak{p}^{3}}{\mathfrak{l}\,^{3}}\Big{(}\frac{X^3-\phi X^2-\xi}{X}\Big{)}\\
&&\mathop{=}^{(\ref{math-p}),(\ref{ties})}
4X^{9}\frac{1000h_{\xi}^{6}(X)f_{\xi}^{3}(X)g_{\xi}^{3}(X)}{X^{12}}\cdot\frac{X^{3}}{8f_{\xi}^{3}(X)},
\end{eqnarray*}
and this equals to the right side of (\ref{math-q}). Thus, we have proved the following. 
\begin{thm}[Triple reduction of the differential system]Let $p(t),q(t),r(t)$ be the solution to the differential system (\ref{sys1}), where $\xi$ is fixed, $-\frac{2+\sqrt{5}}{5}<\xi<\frac{2+\sqrt{5}}{27}$. Let $P=\phi^2 p^2-q^2$, $Q=\phi^2 q^2-r^2$, $R=\phi^2 r^2-p^2$. Let $(X,Y)=(\Upsilon,\Upsilon')$, where $\Upsilon=\Upsilon_{\xi}=\frac{P^3-\phi P^2-\xi}{P}$, or $Q,R$ in the place of $P$. Further, let  
\begin{eqnarray*}
\mathfrak{l}_{\xi}(X)&=&14\phi X-22\xi-4\phi^3,\\
\mathfrak{t}_{\xi}(X)&=&20X^3+5\phi^2 X^2-90\phi\xi X-135\xi^2-20\phi^3\xi.
\end{eqnarray*}
Then a pair of functions $(X,Y)$ parametrizes the algebraic curve
\begin{eqnarray*}
\big{(}Y^2\phi^{5}+\mathfrak{l}_{\xi}(X)\mathfrak{t}_{\xi}(X)\big{)}^2=4\mathfrak{t}_{\xi}^{3}(X).
\label{simpli-i}
\end{eqnarray*}
For $\xi$ in the range we consider, $\xi\neq 0$, this is a curve of arithmetic genus $3$.
\label{thm-red}
\end{thm}
We indeed thus get a triple reduction. Indeed, for $\xi\neq 0$, the triple $(p(t),q(t),r(t))$ parametrizes algebraic curve of arithmetic genus $25$. In Proposition \ref{prop2}, the curve parametrized by $(P,P')$ is of genus $7$. One step further, in Theorem \ref{thm-red} this curve is transformed into a curve of genus $3$. And third, if $\Upsilon(t)$ is any such function, then all three functions $P,Q,R$ are given as three distinct roots of the cubic equation $X^{3}-\phi X^2-\Upsilon(t)X-\xi=0$. Thus, the superflow $\phi_{\mathbb{I}}$ and three functions $p,q,r$ are described in terms of the unique function $\Upsilon$. In case of the octahedral superflow $\phi_{\mathbb{O}}$, this unique function turned out to be the Weierstrass elliptic function. \\

In particular, for $\xi=-\frac{\phi^3}{6}$ we get that $(\Upsilon_{\xi},\Upsilon'_{\xi})$ parametrizes the curve\small
 \begin{eqnarray*}
\Big{(}Y^2\phi^{5}+(14\phi X-\frac{1}{3}\phi^3)
\big{(}20X^3+5\phi^2X^2+15\phi^4 X-\frac{5}{12}\phi^6\big{)}\Big{)}^2=\big{(}20X^3+5\phi^2 X^2+15\phi^4 X-\frac{5}{12}\phi^6\big{)}^{3}.
\end{eqnarray*}
\normalsize
So,  $(\Delta,\Delta')=(\phi^{-2}\Upsilon,\phi^{-2}\Upsilon')$  parametrizes the curve over $\mathbb{Q}$
 \begin{eqnarray*}
4\Big{(}36Y^2+5(42 X-1)
\big{(}48X^3+12X^2+36X-1\big{)}\Big{)}^2=375\big{(}48X^3+12X^2+36 X-1\big{)}^{3}.
\label{special-xi}
\end{eqnarray*}
of arithmetic genus $3$. This is the topic of Section \ref{sec-non-sing}.
\section{A singular case $y=\phi x$}
\label{sec-sing}
Now, we will integrate the vector field in a singular case $y=\phi x$; this is one of the cases corresponding to $\xi=0$. In this case, the vector field reduces to \scriptsize
\begin{eqnarray*}
\frac{2\sqrt{5}(\phi^4x^2-z^2)(\phi^{-2}x^2-z^2)xz}{x^2(1+\phi^2)+z^2}\bl\frac{2\sqrt{5}\phi(\phi^4x^2-z^2)(\phi^{-2}x^2-z^2)xz}
{x^2(1+\phi^2)+z^2}\bl
\frac{-10\phi(\phi^4x^2-z^2)(\phi^{-2}x^2-z^2)x^2}{x^2(1+\phi^2)+z^2}.
\end{eqnarray*}\normalsize
So, the differential system (\ref{sys1}) reduces to
 \begin{eqnarray}
\left\{\begin{array}{l}
p'=-2\sqrt{5}(\phi^4p^2-r^2)(\phi^{-2}p^2-r^2)pr,\\
r'=10\phi(\phi^4p^2-r^2)(\phi^{-2}p^2-r^2)p^2,\\
p^2(1+\phi^2)+r^2=1.
\end{array}
\right.
\label{sys2}
\end{eqnarray}
Of course, the first equality follows from the second and the third. The last one implies $p^2=\frac{1-r^2}{\phi^2+1}$. Thus,
\begin{eqnarray*}
r'&=&10\phi\Big{(}\phi^4\frac{1-r^2}{\phi^2+1}-r^2\Big{)}
\Big{(}\phi^{-2}\frac{1-r^2}{\phi^2+1}-r^2\Big{)}\frac{1-r^2}{\phi^2+1}\Rightarrow\\
r'&=&\frac{2}{\sqrt{5}}(\phi^{2}-4r^2)(\phi^{-2}-4r^2)(1-r^2)=
\frac{2}{\sqrt{5}}(1-12r^2+16r^4)(1-r^2).
\end{eqnarray*}
Integrating, for $r=r(t)$, we obtain:
\begin{eqnarray*}
\frac{(4r^2-2r-1)^2(r+1)}{(4r^2+2r-1)^2(r-1)}=e^{4\sqrt{5}t}.
\end{eqnarray*}
Let $W=e^{4\sqrt{5}t}$. Then the equation for $r$ can be rewritten as
\begin{eqnarray*}
16r^5-20r^3+5r=\frac{W+1}{W-1}:=T,
\end{eqnarray*} 
where $T=\coth(2\sqrt{5}t)$. This can be rewritten as
\begin{eqnarray*}
\prod\limits_{j=0}^{4}\Big{(}r-\sin\Big{(}\frac{2\pi j}{5}\Big{)}\Big{)}=\frac{T}{16}.
\end{eqnarray*}
This identity hints us to use the substitution $r=\frac{m-m^{-1}}{2i}$ for a new algebraic (in $T$) function $m$. This gives
\begin{eqnarray*}
m^5-m^{-5}=\frac{i T}{8}.
\end{eqnarray*} 
This gives 
\begin{eqnarray*}
m^{5}=\frac{\sqrt{256-T^2}+iT}{16}.
\end{eqnarray*}
Let us return to the function $p$ and $r$. Since $p^2=\frac{1-r^2}{1+\phi^2}$, in terms of $m$, this rewrites as
\begin{eqnarray*}
p&=&\frac{m+m^{-1}}{2\sqrt{\phi^2+1}},\\
r&=&\frac{m-m^{-1}}{2i}.
\end{eqnarray*}
This, finally, gives the value
\begin{eqnarray*}
p&=&\frac{1}{\sqrt{10+2\sqrt{5}}}\Bigg{(}\sqrt[5]{\frac{\sqrt{256-T^2}+iT}{16}}+\sqrt[5]{\frac{\sqrt{256-T^2}-iT}{16}}\Bigg{)},\\
r&=&\frac{1}{2i}\Bigg{(}\sqrt[5]{\frac{\sqrt{256-T^2}+iT}{16}}-\sqrt[5]{\frac{\sqrt{256-T^2}-iT}{16}}\Bigg{)},
\end{eqnarray*}
where $T=\coth(2\sqrt{5}t).$
(Compare this to formula (41) in \cite{alkauskas-super1}).
Now,
\begin{eqnarray*}
\phi(x,\phi x,z)=V(x,\phi x,z)\bl V(\phi x,z,x)\bl V(z,x,\phi x).
\end{eqnarray*}
According to a general method to integrate PDE (\ref{pde}),
\begin{eqnarray*}
V(p(s)\v,q(s)\v,r(s)\v)=p(s-\v)\v,\\
V(q(s)\v,r(s)\v,p(s)\v)=q(s-\v)\v,\\
V(r(s)\v,p(s)\v,q(s)\v)=r(s-\v)\v.
\end{eqnarray*} 
In a singular case we are discussing, $x=p(s)\v$, $\phi x=q(s)\v$, $z=r(s)\v$. We have $\v^2=(1+\phi^2)x^2+z^2$.

\section{Generic case for icosahedral superflow, $\xi=-\frac{\phi^{3}}{6}$}
\label{sec-non-sing}
\section{Reducible superflow $\phi_{\mathbb{A}_{4}}$}
\label{sec-anti}


\end{document}